\documentclass{article}
\usepackage[utf8]{inputenc}
\usepackage[a4paper,margin=2.5cm,bmargin=2.5cm]{geometry}
\usepackage{url}

\usepackage[pdftex]{graphicx}
\graphicspath{{../pdf/}{../jpeg/}}
\DeclareGraphicsExtensions{.pdf,.jpeg,.png}

\usepackage{algorithmicx,algpseudocode}
\usepackage[subfigure]{tocloft}
\usepackage{subfigure}
\usepackage{algorithm,algpseudocode}
\usepackage{comment}
\usepackage{amssymb}
\usepackage{amsmath}
\usepackage[english]{babel}
\usepackage{float}
\usepackage{accents}
\usepackage{bm}

\newcommand{\bb}[1]{\mathbf{#1}}
\newcommand{\bs}[1]{\boldsymbol{#1}}
\newcommand{\x}{\bb{x}}

\newcommand{\R}{\mathbb{R}}
\newcommand{\nor}{\bb{n}}
\newcommand\norm[1]{\left\lVert#1\right\rVert}

\newcommand{\Ht}{L^{2}(0,T;H^1(\Omega_{un}))}

\begin{document}

    \title{Fast active thermal cloaking through PDE-constrained optimization and reduced-order modeling}
  \author{Carlo Sinigaglia \thanks{Carlo Sinigaglia is a PhD Candidate at Politecnico di Milano, Department of Mechanical Engineering, Milano 20133, Italy (e-mail: carlo.sinigaglia@polimi.it), Corresponding author.} \and
  Davide E. Quadrelli \thanks{Davide E. Quadrelli is a PhD Candidate at Politecnico di Milano, Department of Mechanical Engineering, Milano 20133, Italy (e-mail: davidee.quadrelli@polimi.it).} \and
    Andrea Manzoni\thanks{Prof. Andrea Manzoni is Associate Professor at Politecnico di Milano, MOX - Department of Mathematics, Milano 20133, Italy (e-mail: andrea1.manzoni@polimi.it).} \and
    Francesco Braghin\thanks{Prof. Francesco Braghin is Full Professor at Politecnico di Milano, Department of Mechanical Engineering, Milano 20133, Italy (e-mail: francesco.braghin@polimi.it). }
  }  




\maketitle

\begin{abstract}
In this paper we show how to efficiently achieve thermal cloaking from a computational standpoint in several virtual scenarios by controlling a distribution of active heat sources. We frame this problem in the setting of PDE-constrained optimization, where the reference field is the solution of the time-dependent heat equation in the absence of the object to cloak. The optimal control problem then aims at actuating the space-time control field so that the thermal field outside the obstacle is indistinguishable from the reference field. In particular, we consider multiple scenarios where material's thermal diffusivity, source intensity and obstacle's temperature are allowed to vary within a user-defined range. To tackle the thermal cloaking problem in a rapid and reliable way, we rely on a parametrized reduced order model built through the reduced basis method, thus entailing huge computational speedups compared to high-fidelity, full-order model exploiting the finite element method while dealing both with complex target shapes and disconnected control domains.  
\end{abstract}

\section{Introduction}
\label{intro}
The possibility of forging a ring or using a cloak to disappear and not be seen by others has stimulated for centuries the imagination of human beings. In practice, an observer that is trying to locate something would not be capable to do it, if the field surrounding the target is locally indistinguishable from that observed without the object itself. The design of a cloaking device implies, by consequence, the capability of covering the target with a special coat, whose properties make the inclusion neutral to the probing field. A solution for such an inverse engineering problem, for long considered impossible to find, has firstly been provided by the development of transformation theories \cite{pendry2006controlling, leonhardt2006optical}. These methods rely on reinterpreting the metric coefficients appearing in the transformed governing equations as space varying material properties, using anisotropy to distort the field inside the cloak, while leaving it unchanged outside.  Originally inspired by the invariance of Maxwell's equation, this theory has spread to all the domains of physics where the governing Partial Differential Equations (PDEs) retain their form under curvilinear coordinates transformations \cite{kadic2015experiments, guenneau2010colours}. Cloaking is thus now an established and active field of research not only in optics, but also in other wave propagation phenomena like  acoustics \cite{cummer2007one,norris2008acoustic,chen2017broadband}, elastodynamics \cite{norris2011elastic,norris2012hyperelastic}, surface water waves \cite{farhat2008broadband} and matter waves \cite{zhang2008cloaking, greenleaf2008approximate}.

Exploiting the close similarity in the mathematical structure of governing equations, the concept of transformation-based cloaks has been extended also to non-wave phenomena, like conduction \cite{schittny2013experiments}, diffusion \cite{guenneau2013fick, schittny2014invisibility} and heat flow \cite{leonhardt2013cloaking}. In hydrodynamics \cite{lamb1924hydrodynamics}, the invariance of the Laplacian operator has long been used to build conformal mappings and solve problems of incompressible flows around complicated shapes exploiting the knowledge of the solution for a simple reference domain. In this context, Transformation Thermodynamics (TT) \cite{guenneau2012transformation} has allowed to achieve experimental evidence of cloaking for thermal fluxes \cite{schittny2013experiments}, based on the implementation of a space varying anisotropic conductivity with alternating layers of indented copper and polydimethylsiloxane (PDMS). A similar approach is used in \cite{ma2013transient}, where the required degrees of freedom for implementing the physical parameter distribution are obtained by adopting five different material ingredients. The weak point of this passive strategy can be thus identified in the high complexity of the resulting cloak to be fabricated: it has been shown that improving the performance of the cloak designed in \cite{schittny2013experiments} implies increasing the number of alternating material's layer of three orders of magnitude \cite{petiteau2014spectral}. Moreover, the shapes of the target for which solutions are readily available are limited to simple geometries. Efforts have thus been put on finding alternative strategies to overcome these weaknesses. As an example, Ji et al \cite{Ji2019} have recently compared TT to the Neutral Inclusion method \cite{kerner1956elastic}, which allows for cloaking with a single layer of isotropic material. 

Among alternative methods, another possible option is represented by the so-called "active" cloaking: firstly developed in the context of the Helmholtz equation \cite{miller2006perfect, vasquez2009active, vasquez2009broadband, vasquez2011exterior}, where sensors and sources are used in a similar fashion as in classical noise cancellation \cite{lin2021active}, it has become an attractive solution for the thermal case as well \cite{nguyen2015active,cassier2021active,Huang2019}. Indeed, thermoelectric materials can be used to design heat sources and sinks that exert a suitable manipulation of the thermal field in a localized region of the domain (usually around the target to cloak) to achieve the cloaking objective. In \cite{nguyen2015active}, an air hole in a steel plate is thermally cloaked in the steady-state regime by using Peltier elements, while heat flows from the left to the right side of a steel plate. 

More recently, in \cite{cassier2021active} an active thermal cloaking technique is developed exploiting Green representation formulas, and is able to cloak the target in the transient regime. In this case, the unperturbed transient evolution of the thermal field due to a heat source in an unbounded domain is described through boundary integrals defined on the active cloak region; this in turns allows to generate the space-time modulated actuation needed to hide an obstacle that would otherwise perturb the thermal background field. The same approach can be used for mimicking, that is, the problem of making the thermal trace of an object appear as the one of a completely different target. However, the approach introduced in \cite{cassier2021active} has the drawback of dealing with complex boundary integrals and does not take into account the dimensionality of the thermoelectric actuators. Furthermore, it requires a closed surface for the active cloaking to work. \\

In this paper, we take instead another route to obtain thermal cloaking, that is, we tackle the design phase by reformulating the problem following an optimization strategy. We formulate an infinite-dimensional optimization problem whose objective is aimed at tracking the unperturbed thermal field generated by a heat source. Similar ideas from PDE-constrained optimization were exploited in \cite{chen2021optimal} to derive material properties for the design of passive acoustic cloaks while in \cite{cominelli2021design} an optimization strategy is derived to design cloaks implementable with sonic metamaterials. We frame this objective in the context of optimal control of PDEs where the state equation describes the temperature field in presence of an obstacle with constant temperature (i.e., the heat diffusion equation), and the active cloaking sources, playing the role of distributed control functions, are defined in an annular region around the obstacle to cloak. With respect to other existing techniques, this allows us to  show that connectedness of the control domain is not necessary to achieve the cloaking objective, thus setting the ground towards more realistic applications. Furthermore, the flexibility of the optimization framework allows to cloak complex objects both in the steady-state and the transient regime. We first consider the steady-state problem and then proceed to the transient case. After deriving the optimality system of equations to be fulfilled by the optimal state and actuation fields, we derive its numerical approximation relying on a high-fidelity, full-order model consisting of a finite element discretization in space.
Indeed, as most other active cloaking strategies, our approach presents the significant drawback that the probing reference field needs to be known ahead of time. However, due to the computational efficiency and real-time potential of our strategy, an estimation problem could be implemented with similar optimization methods thus overcoming this important issue by reconstructing the probing field with real-time measurements.


 At the discrete level, we propose a fast and efficient way to solve in one-shot the steady-state problem and present the computational issues in treating the transient problem. Then, we develop a Reduced Order Model (ROM) to reduce the computational burden of the Optimal Control Problem (OCP) for a different set of scenario parameters, thus dealing with the model order reduction of a class of \emph{parametrized} OCPs (see e.g., \cite{red_book}). In particular, we make use of the Proper Orthogonal Decomposition (POD) technique to simultaneously reducing the reference and the optimal control problem while considering a set of scenario parameters that affect the solution both of the reference and of the optimal control. In particular, we are interested in studying parametric OCPs where obstacle's temperature, thermal diffusivity and source intensity are parameterized. The POD reduction methods is particularly well suited for parabolic problems (see e.g., \cite{red_book,heat_POD,azaiez2016recursive}) since the diffusion operator filters out higher frequency contributions and a lower dimensional solution can be well resolved with carefully selected reduced basis, if the dependence on parameters is sufficiently smooth \cite{red_book}. In this way, it is possible to recompute the optimal control field for a new set of such parameters in a fraction of a second. The computational speedup achieved enables to deploy the control algorithm in situations where the obstacle temperature or the material properties of the background are not precisely known but they are known to belong to some parameter set.  \\

 The paper is organized as follows. In Section~\ref{OCP}, the cloaking problem is formulated as a PDE-constrained optimization problem and a set of first-order optimality conditions is derived. In Section~\ref{ROM} the reduced-order model for both the reference and the optimal control problem is developed. In Section~\ref{NUM}, several test cases both for the steady-state and the transient problem are presented to show the effectiveness of the proposed method, while some conclusions are finally reported in Section~\ref{conclusion}.

\section{The optimal control problem}
\label{OCP}
In this Section, we formulate the cloaking objective as an OCP in the PDE settings. As a reference unperturbed state, we consider the thermal field generated by a source in a bounded domain $\Omega_{un}$ during the time interval $(0,T)$. The reference thus satisfies the heat equation: \vspace{-0.1cm}
\begin{equation}
\label{heat_basic}
   \displaystyle \rho c \frac{\partial z(\x,t)}{\partial t} - k \Delta z(\x,t) = \tilde{s}(\x) \quad \text{in} \quad \Omega_{un}\times (0,T) \vspace{-0.1cm}
\end{equation}
where the temperature $z$ is measured in Kelvin, $\x \in \R^2$ is the spatial coordinate in meters, $t$ indicates time in seconds, $c$ is the specific heat ($\text{J}\text{K}^{-1}\text{kg}^{-1}$), $\rho$ the mass density ($\text{kg}\,\text{m}^{-3}$) and $k$ is the thermal conductivity ($\text{W}\,\text{m}^{-1}\,\text{K}^{-1}$) of the material. It is useful to define the thermal diffusivity $\mu = \frac{k}{\rho c_p}$ ($\text{m}^2\text{s}^{-1}$) and rewrite equation (\ref{heat_basic}) equivalently as: \vspace{-0.1cm}
\begin{equation}
\label{heat_basic}
\displaystyle     \frac{\partial z}{\partial t} - \mu \Delta z = s  \quad \text{in} \quad \Omega_{un}\times (0,T) \vspace{-0.1cm}
\end{equation}
where $s = \frac{\tilde{s}}{\rho c_p}$ ($\text{K}\text{s}^{-1}$) is the source term.
To solve it, we equip the problem with suitable initial and boundary conditions. For the sake of simplicity, we consider a square computational domain with homogeneous Robin boundary conditions that approximate to first-order an unbounded domain \cite{ABCs}. Hence, the unperturbed field $z \in \Ht$ is the solution of: \vspace{-0.1cm}
\begin{equation}
\label{reference_state}
\begin{cases}
\begin{array}{ll}
\frac{\partial z}{\partial t} -\mu\Delta z    = s &  \textrm{in} \quad \Omega_{un}\times (0,T) 
\phantom{space} \\
-\mu\nabla z \cdot \nor + \alpha \, z = 0 & \textrm{on} \quad \partial \Omega_{un} \times (0,T)  \\
z(\x,0)= 0 & \textrm{in} \quad \Omega_{un}\times \{0\}   \\
\end{array}
\end{cases} \vspace{-0.1cm}
\end{equation}
whereas the steady-state reference is simply obtained by considering the Poisson equation $ -\mu\Delta z    = s$ with the same boundary conditions on $\Gamma_d = \partial \Omega_{un}$ and neglecting the initial conditions. 

\begin{figure}[h!]
\centering
\subfigure{\includegraphics[width=0.7\textwidth]{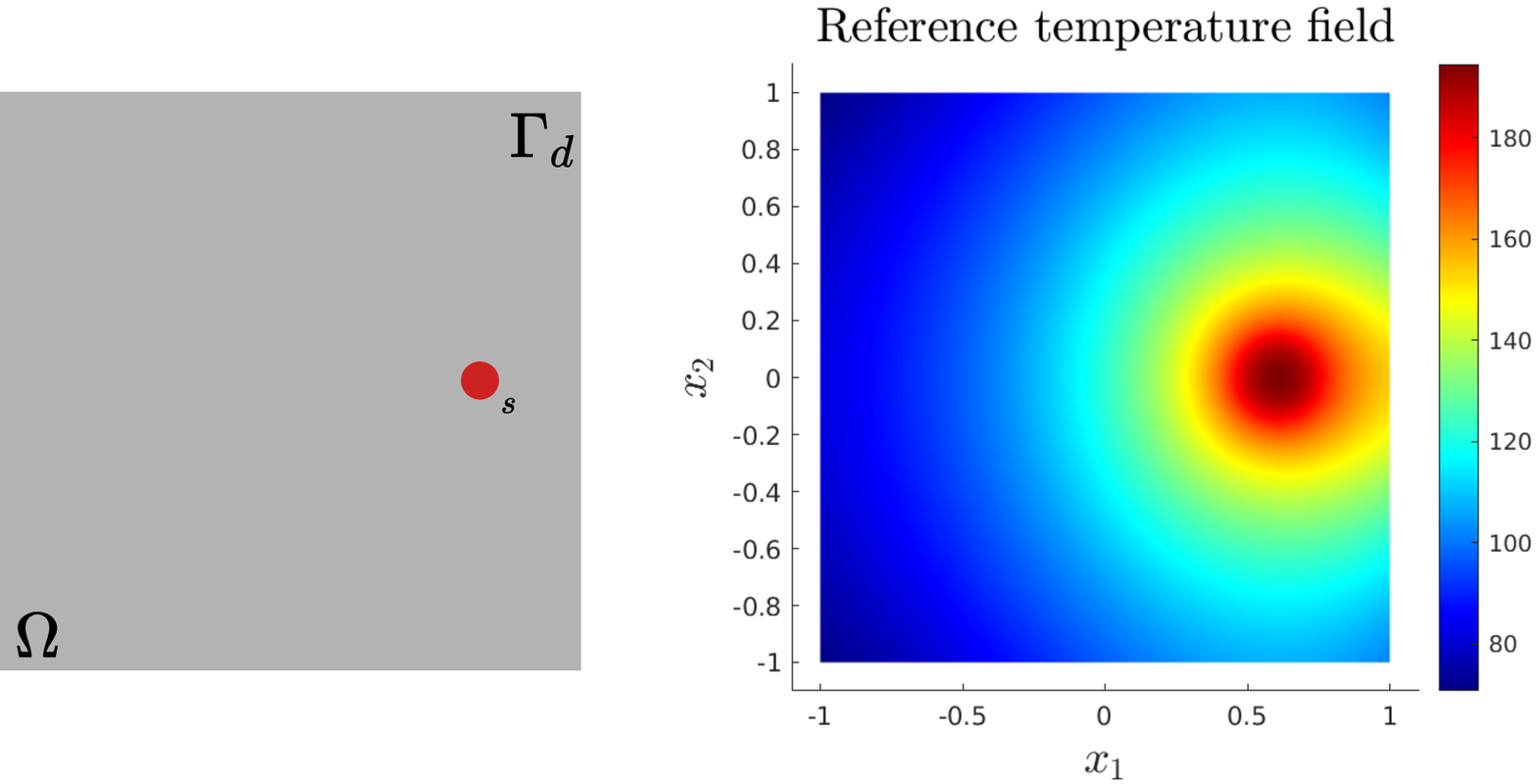}}
\subfigure{\includegraphics[width=0.7\textwidth]{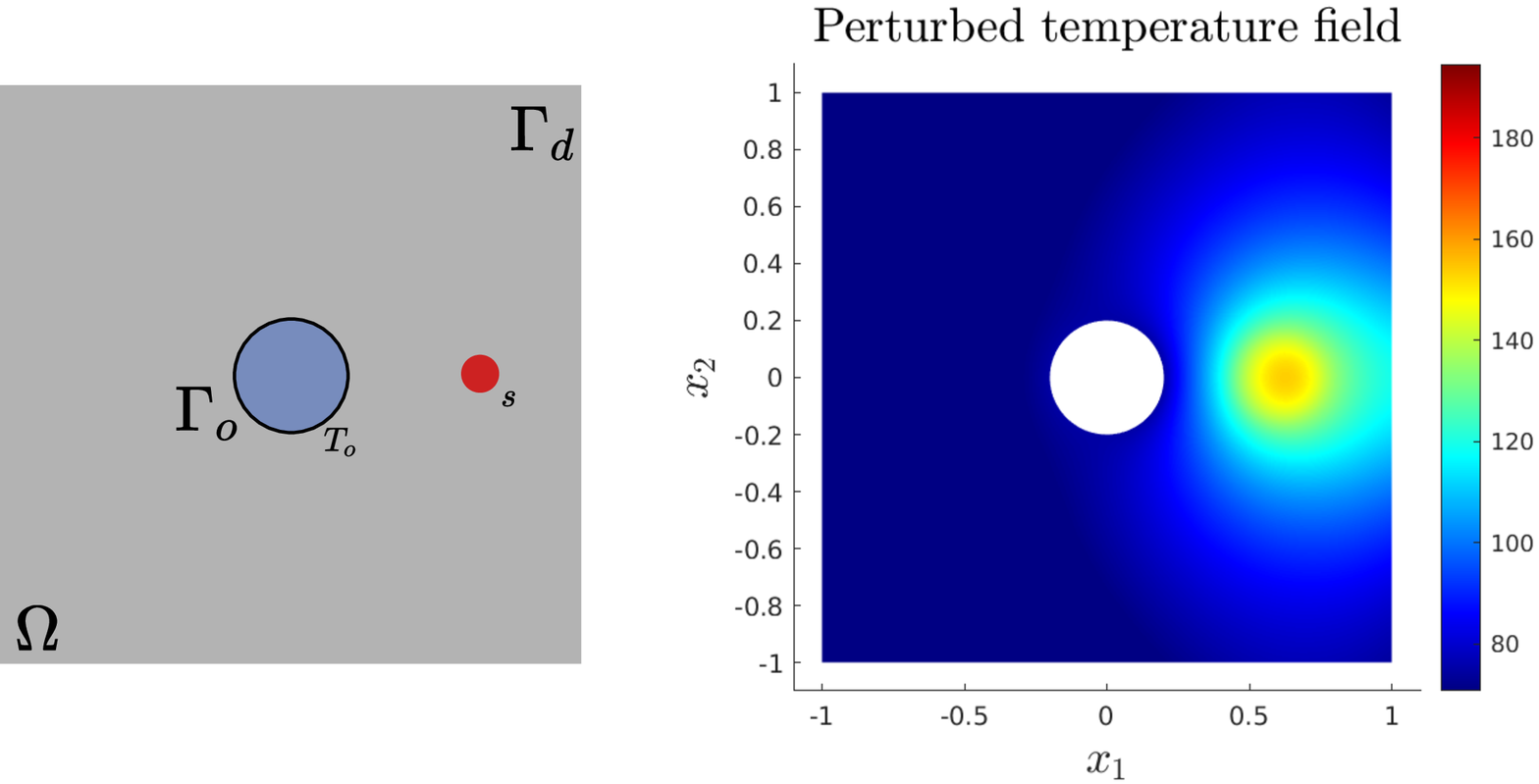}}
\caption{Top Figure: Unperturbed temperature field generated by a heat probing source at the mid-right part of the domain for $\mu = 5$, $\alpha=1$ and $s=10000$. Bottom Figure: Perturbed temperature field generated by a certain obstacle with constant temperature $T_o = 75^{\circ}C$. The obstacle generates a thermal signature that modifies the temperature field in its neighborhood. The parameters $\mu$,$\alpha$ and $s$ are the same as the top figure. }
\label{reference_layout_fig}
\label{perturbed_state}
\end{figure}

Note that in order to simulate the domain's unboundedness, the simple choice $\alpha=1$ allows to obtain good results in practice \cite{ABCs}.  The unperturbed layout together with the steady-state solution obtained through the FEM for selected values of $\mu$,$\alpha$,$s$ is shown in Figure \ref{reference_layout_fig}.

 We now consider the temperature field generated by considering the presence of an obstacle $\Theta \subset \Omega_{un}$ whose temperature is kept constant. The obstacle thus induces an inhomogeneous boundary conditions of Dirichlet type along its boundary $\Gamma_o = \partial \Theta$ that perturbs the background thermal field. We denote as $q \in {L^{2}(0,T;H_{\Gamma_o}^1(\Omega))}$ the temperature field in the presence of the obstacle, that is the solution of: \vspace{-0.1cm}
\begin{equation}
\label{state_problem}
\begin{cases}
\begin{array}{ll}
\frac{\partial q}{\partial t} -\mu\Delta q    = s + u &  \textrm{in} \quad \Omega \times (0,T) 
\phantom{space} \\
-\mu\nabla q \cdot \nor + q = 0 & \textrm{on} \quad \Gamma_d \times (0,T) \\
q = T_o & \textrm{on} \quad \Gamma_o \times (0,T) \\
q(\x,0)= 0 & \textrm{in} \quad \Omega \times \{0\}.  \\
\end{array}
\end{cases}
\end{equation}
where $\Omega = \Omega_{un} \setminus \Theta $ and $\partial \Omega = \Gamma_o \cup \Gamma_d $.
A layout of the state problem and the effect of the Dirichlet obstacle is shown in Figure \ref{perturbed_state} where an obstacle's temperature of $T_o = 75^{\circ}C$ is considered. Note that the farther the temperature $T_o$ from the temperature in the corresponding region of the unperturbed field, the more evident the effect of the perturbation.

\begin{figure*}[hbt!]
\centering
\includegraphics[width=0.7\textwidth]{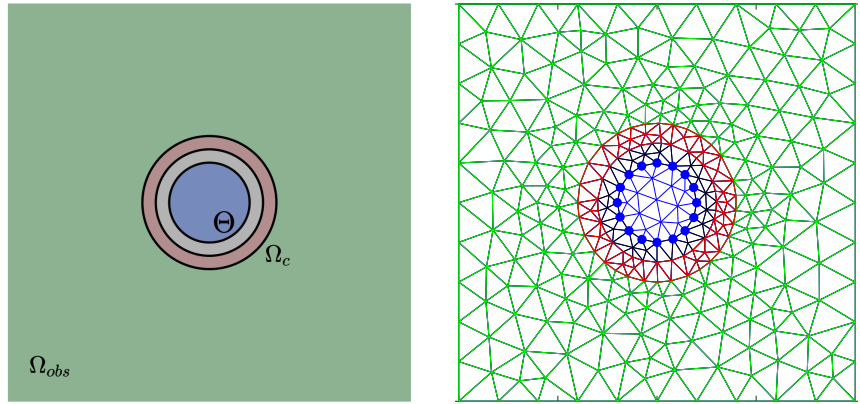}
	\caption{Conceptual domains definition (left plot) and illustrative coarse computational mesh (right plot). The unperturbed problem is defined on the square $\Omega_{un}$ while the domain of the control problem is $\Omega = \Omega_{un} \setminus \Theta$. The control domain $\Omega_c$ and the observation domain $\Omega_{obs}$ are shown in red and green, respectively. The obstacle domain $\Theta$ is colored in blue. The Dirichlet boundary condition induced by $\Theta$ is represented by the blue points on the right. }
	\label{control_layout}
\end{figure*}

We now turn to the setup of the Optimal Control Problem (OCP), referring to Figure \ref{control_layout} for its illustration. We consider an annular region $\Omega_c \subset \Omega$ in which a space-time varying control source $u(\x,t)$ is defined.
The cloaking objective can be rephrased as finding $u \in L^2(0,T,H^1(\Omega_c))$ so that the state $q$ is as close as possible to the reference thermal field $z$ onto an observation region $\Omega_{obs} \subset \Omega$, exterior to the object. In our case the control region is defined as a circular annulus with thickness $r_c$ surrounding the target to cloak. Note that the flexibility allowed by the PDE-constrained formulation allows the control region to be arbitrarily shaped in the domain thus bridging the gap towards realistic applications.

Therefore, the cloaking problem can be recast as the following OCP: \vspace{-0.1cm}
\begin{equation}
\label{ocp_formulation}
\begin{aligned}
&\min_{u,q} \quad  J(q,u) =  \frac{1}{2} \int_{0}^T \int_{\Omega_{obs}}  (q(\x,t)-z(\x,t))^2 \, d \Omega \, dt    + \frac{1}{2} \int_{0}^T \int_{\Omega_c} \Big(\beta\,u(\x,t)^2  \, + \beta_g\,||\nabla u(\x,t)||^2 \Big)     d \Omega d t                  \\
&  \mbox{subject to} \\
& 
\begin{cases}
\begin{array}{ll}
\frac{\partial q}{\partial t} -\mu\Delta q    = u + s  &  \textrm{in} \quad \Omega \times (0,T) 
\phantom{space} \\
-\mu\nabla q \cdot \nor + q = 0 & \textrm{on} \quad \Gamma_d \times (0,T) \\
q = T_o            & \textrm{on} \quad \Gamma_o \times (0,T) \\
q(\x,0)=0 & \textrm{in} \quad \Omega \times \{0\},  \\
\end{array}
\end{cases}
\end{aligned}
\end{equation}
where $z$ is the solution of the reference problem (\ref{reference_state}). The OCP (\ref{ocp_formulation}) consists of the minimization of a quadratic cost functional subject to a linear constraint in the form of the parabolic heat equation.

The theory of linear-quadratic OCPs in the PDE setting is well developed (see e.g. \cite{fredi,salsa,quart}) and the existence and uniqueness of a solution can be proved. Reasoning in a similar way, we can state the steady-state OCP as: \vspace{-0.1cm}
\begin{equation}
\label{SS_ocp_formulation}
\begin{aligned}
&\min_{u,q} \quad  J_{ss}(q^{ss},u^{ss})  =  \frac{1}{2}  \int_{\Omega_{obs}}  (q^{ss}(\x)-z^{ss}(\x))^2 \, d \Omega    + \frac{1}{2}  \int_{\Omega_c} \Big(\beta\,u^{ss}(\x)^2  + \beta_g\,||\nabla u^{ss}(\x)||^2\Big)  \,    d \Omega                 \\
&  \mbox{subject to} \\
& 
\begin{cases}
\begin{array}{ll}
-\mu\Delta q^{ss}    = u^{ss} + s  &  \textrm{in} \quad \Omega \\
-\mu\nabla q^{ss} \cdot \nor + q^{ss} = 0 & \textrm{on} \quad \Gamma_d \\ 
q^{ss} = T_o            & \textrm{on} \quad \Gamma_o.  \\
\end{array}
\end{cases}
\end{aligned}
\end{equation}

We now proceed by recovering a system of first-order necessary and sufficient optimality conditions for the transient problem, exploiting the Lagrange method \cite{fredi,herzog}.  The Lagrangian functional $\mathcal{L}: \mathcal{V} \times \mathcal{U} \times \mathcal{W}^{*} \mapsto \R$ is defined as \vspace{-0.1cm}
\begin{equation*}
    \mathcal{L}(q,u,p) = J(q,u) + \langle p , G(q,u) \rangle_{\mathcal{W}^*,\mathcal{W}} \vspace{-0.1cm}
\end{equation*}
where $\mathcal{V}=H^{1}(0,T;H_{\Gamma_o}^{1}(\Omega),H_{\Gamma_o}^{1}(\Omega)^*)$, $\mathcal{W}=L^2(0,T;H_{\Gamma_o}^1(\Omega)^{*})$ and $G$ is the abstract PDE constraint, so that state, control and adjoint variables are considered independently. Therefore, the set of first-order necessary optimality conditions consists of imposing that the Gateaux derivative of the Lagrangian with respect to the triple $(q,u,p)$ along an arbitrary variation $(\psi,h,\phi)$ is equal to zero. 
In our case, the Lagrangian takes the explicit form \vspace{-0.1cm}
\begin{equation}
\label{lag_def}
    \mathcal{L} = J(q,u) + \int_{\Omega}\int_0^T \left( -\frac{\partial q}{\partial t} +\mu\Delta q +  u + s \right) p  \,\, d \Omega d t , \vspace{-0.1cm}
\end{equation}
and the resulting optimality system in strong form consists of the state dynamics (\ref{state_problem}), the adjoint dynamics: \vspace{-0.1cm}
\begin{equation}
\label{adjoint_problem}
\begin{cases}
\begin{array}{ll}
-\frac{\partial p}{\partial t} -\mu\Delta p    = (q-z)\chi_{obs} &  \textrm{in} \quad \Omega \times (0,T) 
\phantom{space} \\
-\mu\nabla p \cdot \nor + p = 0 & \textrm{on} \quad \Gamma_d \times (0,T) \\
p = 0 & \textrm{on} \quad \Gamma_o \times (0,T) \\
p(\x,T)= 0 & \textrm{in} \quad \Omega \times \{0\}  \\
\end{array}
\end{cases}
\end{equation}

\noindent and the optimality condition (Euler equation) :
\begin{equation*}
\begin{array}{ll}
-\beta_g \Delta u + \beta u + p    = 0 &  \textrm{in} \quad \Omega_c \times (0,T) .
\end{array}
\end{equation*}

\noindent In the adjoint problem (\ref{adjoint_problem}), $\chi_{obs}$ denotes the characteristic function of the observation domain $\Omega_{obs}$ that is equal to one on that domain and zero elsewhere. Furthermore, we remark that the adjoint variable $p$ satisfies the backwards heat equation which is ill-posed when the time datum is given as an initial condition. In our case however, the problem is well-posed since a terminal condition is prescribed (see e.g. \cite{fredi}, Chapter 3).

\noindent Since this kind of OCPs is rather standard, we skip the details about the derivation of the optimality system and refer the reader to \cite{fredi}, Chapter 3. On the other hand, the system of optimality conditions for the steady-state problem is obtained simply by deleting the time derivatives in the state and the adjoint problem. Note that the transient problem is cast along a finite time horizon $(0,T)$ without any terminal cost, that is, we did not put any weight on the final state $q(\x,T)$ and control $u(\x,T)$ achieved. As a consequence, the final adjoint condition is homogeneous $p(\x,T)=0$. However, it is reasonable to require that the optimal solution of the transient OCP (\ref{ocp_formulation}) converges to the optimal solution of the steady-state OCP (\ref{SS_ocp_formulation}). This condition is enforced by replacing the final adjoint condition with $p(\x,T)= {p}^{ss}(\x)$, where ${p}^{ss}$ is the solution of the steady-state OCP. 
Finally, to handle several possible scenarios at the same time, we consider as input parameters the diffusivity constant $\mu$, the obstacle's temperature $T_o$ and the source intensity $I$, thus ending up with two \emph{parametric} OCPs for the steady-state and the transient heat transfer problems coupled with the reference dynamics $z$ that in turn depends on the input parameters as well. Note that the new final adjoint condition generates a one-way coupling between the steady-state and the transient problems.

\subsection{Numerical discretization of the OCP}

For the numerical approximation of the reference problem and the OCPs, we employ the Finite Element Method (FEM). We introduce a triangulation of the domain $\Omega_{un}$ made by triangular elements and select piecewise linear, globally continuous ansatz functions $\phi_i$ ($\mathbb{P}_1$ finite elements) for the space approximation of state and adjoint variables in $\Omega$, while the basis are restricted to $\Omega_c$ for the control approximation. 

To generate a suitable triangulation on $\Omega$, we restrict the previous mesh to those elements belonging to $\Omega_{un} \setminus \Theta$, so that the nodes that belong to the obstacle are masked and a Dirichlet boundary is introduced instead. An illustrative example of the computational mesh and the resulting domain definition is shown in Figure \ref{control_layout}. We define by $E\in \R^{N_q \times N_z}$ the matrix that restricts the reference nodes to the OCP mesh nodes, where $N_q$ and $N_z$ are the dimension of the state and reference variables, respectively. In this way, a matrix $A\in \R^{N_z \times N_z}$ assembled on the reference mesh can be restricted to the OCP mesh through the projection $\tilde{A} = E A E^{\top}$.
In the following, we indicate with a tilde the matrices restricted to the computational domain of the control problem, e.g, $\tilde{A} = EAE^{\top}$.
The FEM approximation of the reference dynamics is therefore \vspace{-0.1cm}
\begin{equation}
\label{reference_FEM}
    \begin{cases}
    M \dot{\bb{z}} + A\bb{z} = \bb{F}  \, , \, \qquad t \in (0,T) \\
    \bb{z}(0)                = \bb{0}  \, , \\ 
    \end{cases}
\end{equation}
where $M_{ij} = \int_{\Omega_{un}} \phi_i \phi_j d\Omega$, $A_{ij} = \int_{\Omega_{un}} \mu \nabla \phi_i \cdot \nabla \phi_j d\Omega + \int_{\partial \Omega_{un}} \phi_i \phi_j d\Gamma$ and $F_i = \int_{\Omega_{un}} s \phi_i d\Omega$ for $i,j = 1\, \ldots, N_z$. Regarding the OCP discretization, we follow an Optimize-then-Discretize strategy (see e.g. \cite{herzog}) that consists of discretizing the optimality conditions obtained at the continuous level. Hence, the discretization of the optimality system gives: \vspace{-0.1cm}
\begin{equation}
    \begin{cases}
         EME^{\top}  \dot{\bb{q}} + EAE^{\top}\bb{q} = \bb{F}_o + E\bb{F} + B\bb{u}  \, , \, & \quad t \in (0,T)  \\
     \bb{q}(0) = \bb{0}  &  \\
    -EME^{\top}  \dot{\bb{p}} + EAE^{\top}\bb{p} = M_{obs}(\bb{q}-E\bb{z}) \, , \, & \quad  t \in (0,T)  \\ 
    \bb{p}(T)  = \bb{p}^{ss} & \\
    (\beta M_u+\beta_g A_u) \bb{u} + B^{\top} \bb{p}          = \bb{0} \, , \, & \quad  t \in (0,T) 
    \end{cases}
    \label{opt_discrete}
\end{equation}
where $\bb{F}_o$ is the contribution arising from the Dirichlet data, $B\in \R^{N_q \times N_u}$ is the control matrix with elements $B_{ik}= \int_{\Omega} \phi_i \phi_k d\Omega$, while $M_{obs}$ is the observation domain mass matrix, $M_u$ and $A_u$ the control mass and diffusion matrices. The role of $A_u$ and $\beta_g$ will be made clear later on. 

On the other hand, the algebraic optimality system for the steady-state problem \eqref{SS_ocp_formulation} -- whose variables are denoted with a ``ss" superscript -- reads
\begin{equation}
\label{steady_state_FEM}
    \begin{cases}
     A\bb{z}^{ss} = \bb{F}  \\
     EAE^{\top}\bb{q}^{ss} = \bb{F}_o + E\bb{F} + B\bb{u}^{ss}  \\
     EAE^{\top}\bb{p}^{ss} = M_{obs}(\bb{q}^{ss}-E\bb{z}^{ss}) \\ 
    (\beta M_u+\beta_g A_u) \bb{u}^{ss} + B^{\top} \bb{p}^{ss}          = \bb{0}.
    \end{cases}
\end{equation}
Note that  system (\ref{steady_state_FEM}) can be solved with a ``one-shot" method by grouping the unknowns as: \vspace{-0.1cm}
\begin{equation}
\label{big_ss}
\underbrace{\begin{bmatrix}
A & \phantom{} & \phantom{} & \phantom{} \\
\phantom{} & EAE^{\top} & \phantom{} & -B \\
M_{obs} E      & -M_{obs}       & EAE^{\top} & \phantom{} \\
\phantom{} & \phantom{} & B^{\top}   & \beta M_u +\beta_g A_u
\end{bmatrix}}_{\mathcal{K}}
\underbrace{\begin{bmatrix}
\bb{z}^{ss} \\
\bb{q}^{ss} \\
\bb{p}^{ss} \\
\bb{u}^{ss} 
\end{bmatrix}}_{\bb{y}}
=
\underbrace{\begin{bmatrix}
\bb{F} \\
\bb{F}_o + E\bb{F} \\
\bb{0} \\
\bb{0} 
\end{bmatrix}}_{\mathcal{F}} \vspace{-0.1cm} 
\end{equation}
so that the reference state, optimal state, adjoint and control variables are obtained from the solution of a single (yet moderately large) linear system with a saddle-point structure.

 Regarding the time discretization of the transient problems, we rely on the Crank-Nicolson method, so that the systems of ODEs involved in (\ref{reference_FEM}) and (\ref{opt_discrete}) are turned into large algebraic systems to be solved at each time step. We note that, being the adjoint equation well-posed backwards in time, there are no issues in the time discretization. We discretize the time interval $(0,T)$ into $N$ subintervals of length $\Delta t = \frac{T}{N}$ and we use the subscript $k=0,\ldots,N$ to indicate a vector computed at a discrete time instance $t_k=k\Delta t$.
Thus, we look for the unknown vectors $\bb{z}_k,\bb{q}_k,\bb{p}_k,\bb{u}_k$ at each instant $t_k$ which solve the discretized version of the optimality system (\ref{opt_discrete}). Although a one-shot implementation is also possible \cite{Strazzullo2020}, we rely on an iterative quasi-Newton method that allows for a finer discretization of the OCP. Indeed, the dimension of the linear system in that case would result $O(N(N_z+2N_q + N_u))$, thus making its handling computationally intensive. The pseudocode is shown in Algorithm (\ref{alg_newton}). 
\begin{algorithm}
  \begin{algorithmic}[1]
  \State $\bb{z}^{ss},\bb{q}^{ss},\bb{p}^{ss},\bb{u}^{ss} \gets \text{Solve ``one-shot" linear system (\ref{big_ss}) } $ \Comment{Solve steady-state OCP}
  \State $\bb{z} \gets \text{Solve reference dynamics} $
  \State $\bb{p}_N \gets \bb{p}^{ss}$ \Comment{Assign final condition to adjoint variable}
  \State $\bb{u}^{(0)}_k \gets \bb{u}^{ss} \quad \text{for} \,\, k = 0 , \ldots , N$ \Comment{Initialize control variable as the steady-state one}
    \For{$i=0:\text{maxIter}$} 
        \State    $\bb{q}^{(i)} \gets \text{Solve state equation with control } \bb{u}^{(i)}$ 
        \State    $\bb{p}^{(i)} \gets \text{Solve modified adjoint equation with control } \bb{u}^{(i)} \text{ and state } \bb{q}^{(i)}$ 
        \State $\nabla J(\bb{u}^{(i)}) \gets (\beta M_u+\beta_g A_u) \bb{u}^{(i)} + B^{\top} \bb{p}^{(i)}$ \Comment{Compute reduced gradient}
        \State $\bb{d}^{(i)} \gets \text{Solve linear system } (\beta M_u+\beta_g A_u) \bb{d}^{(i)} = -\nabla J(\bb{u}^{(i)})  $ \Comment{Quasi-Newton direction}
        \State $\tau \gets \text{ArmijoBacktracking}(J,\bb{d}^{(i)},\bb{u}^{(i)}) $ \Comment{Line search in direction $\bb{d}^{(i)}$}
        \State $\bb{u}^{(i+1)} \gets \bb{u}^{(i)} + \tau \bb{d}^{(i)}  $ \Comment{Update control}
        \If{ $\norm{\nabla J(\bb{u}^{(i)})} < \text{tol}$}
        \State $  \text{return}$
        \EndIf
    \EndFor
  \end{algorithmic}
  \caption{Modified Newton method for parabolic OCPs  }
  \label{alg_newton}
\end{algorithm}
The Quasi-Newton step for the control update at each iteration can be written as:
\begin{equation*}
\begin{aligned}
      \bb{d}_k^{(i)} &= -(\beta M_u+\beta_g A_u)^{-1}\Big((\beta M_u+\beta_g A_u) \bb{u}_k^{(i)} + B^{\top} \bb{p}_k^{(i)}\Big) \\
      &= -\bb{u}_k^{(i)} - (\beta M_u+\beta_g A_u)^{-1} B^{\top} \bb{p}_k^{(i)}
\end{aligned}   
\end{equation*}
for time index  $k=0,\ldots,N$. Then the update reads:
\begin{equation}
\label{fixed_step}
\bb{u}^{(i+1)}_k = (1-\tau)\bb{u}^{(i)}_k -\tau (\beta M_u+\beta_g A_u)^{-1} B^{\top} \bb{p}^{(i)}_k=(1-\tau)\bb{u}^{(i)}_k+\tau \Lambda \bb{p}_k^{(i)} 
\end{equation}
where we have defined the matrix $\Lambda = -(\beta M_u+\beta_g A_u)^{-1} B^{\top}$. We also remark the fact that the update step (\ref{fixed_step}) can be interpreted as a fixed-point iteration and that the preconditioner $(\beta M_u+\beta_g A_u)$ does not correspond to the reduced Hessian, that is the Hessian of the OCP interpreted as a function of the control variable $\bb{u}$ only.

\section{Reduced order model}
\label{ROM}
In this Section, the optimal control formulation of the cloaking problem is considered in a parametrized setting, to take into account a range of possible scenarios of physical interest \cite{Kahlbacher2012,Kunisch2001,Benner2014}, and determine the optimal cloaking strategy for each scenario inexpensively. Indeed, relying on high-fidelity solvers -- such as the ones that exploit the finite element method -- would be infeasible if the OCP must be solved several times, or in a very short amount of time. To this goal, we rely on reduced order modeling techniques for parametrized PDEs, applying the reduced basis (RB) method to the optimality system arising from the  (now, parametrized)  OCP. This represents a further step towards practical applications of the cloaking algorithm, since  a ROM enables a reduction of the computational complexity by orders of magnitude.

The parametrized version of the optimality system \eqref{opt_discrete} and the reference dynamics \eqref{reference_FEM} for the case at hand can be written as: \vspace{-0.1cm}
\begin{equation}
    \begin{cases}
        M \dot{\bb{z}} + A(\bs{\mu})\bb{z} = \bb{F}(\bs{\mu}) \, , \, & t \in (0,T) \\
        \bb{z}(0)                = \bb{0}  & \smallskip \\ 
        \tilde{M}  \dot{\bb{q}} + \tilde{A}(\bs{\mu})\bb{q} = E\bb{F}(\bs{\mu}) + \bb{F}_{o}(\bs{\mu}) + B\bb{u}   \, , \, &t \in (0,T) \\
     \bb{q}(0) = \bb{0} & \smallskip \\
    -\tilde{M}   \dot{\bb{p}} + \tilde{A}(\bs{\mu})\bb{p} = M_{obs}(\bb{q}-E\bb{z})  \, , \, & t \in (0,T) \\ 
    \bb{p}(T)  = \bb{p}^{ss} & \smallskip \\
    (\beta M_u + \beta_g A_u) \bb{u} + B^{\top} \bb{p}          = \bb{0}  \, , \, & t \in (0,T)
    \end{cases} \vspace{-0.1cm}
    \label{eq:param_opt_discrete}
\end{equation}
where $\bs{\mu}=[ \, \mu \,\, I \,\, T_{o} \, ]$ and $\bb{F}_{o}$ is the forcing term arising from the non-zero temperature of the obstacle. In the system above we have highlighted which matrices depend on the parameter vector $\bs{\mu}$. The set of parameters considered are the temperature of the obstacle $T_o$, the thermal diffusivity of the material $\mu$ and the intensity of the heat source $I$. The parameter space, representing the (product) of the range of variation of each parameter, is denoted by $\mathcal{P} \subset \mathbb{R}^p$; in our case $p=3$. Note that  we still treat the control variable as a distributed field, and describe through a set of input parameters a range of virtual scenarios of interest. Similarly, the matrix $A$ and the vector ${\bf F}$ would depend on  $\bs{\mu}=[ \, \mu \,\, I \,\, T_{o} \, ]$, with the same meaning of the input parameters, in the parametrized version of the optimality system \eqref{steady_state_FEM} for the steady OCP.


\subsection{The Reduced Basis Method for linear transient parametrized problems}

We now briefly recall how to set up a ROM through the RB method for a linear transient parametrized problem, as in the case of the state and the adjoint problems appearing in the system \eqref{eq:param_opt_discrete}. For the sake of simplicity, we focus on the case of a $N_y$-dimensional problem of the form  \vspace{-0.1cm}
\begin{equation}
    \begin{cases}
        {M}  \dot{\bb{y}} +  {A}(\bs{\mu})\bb{y} = \bb{F}(\bs{\mu})    \, , \, &t \in (0,T) \\
     \bb{y}(0) = \bb{y}_0 
    \end{cases} \vspace{-0.1cm}
    \label{eq:param_sample_pb}
\end{equation}
and postpone to the following section further details about the coupling between state, adjoint, and control variables. A POD-Galerkin RB method to approximate the solution of the FOM \eqref{eq:param_sample_pb} aims at approximating ${\bb{y}} = {\bb{y}}(t, \bs{\mu})$ as a linear combination of $n \ll N_y$ POD modes, under the form ${\bb{y}}(t, \bs{\mu}) \approx V  {\bb{y}}_n(t, \bs{\mu})$. Here $V \in \mathbb{R}^{N_y \times n}$ denotes a matrix whose columns are orthogonal and represent the $n$ basis vectors spanning the reduced basis subspace the solution is sought in. This reduced basis is built from a set of $ n \times n_s$ snapshots $\{ {\bb{y}}(t_k, \bs{\mu}^i), \ k=1,\ldots,N, \, i=1,\ldots, n_s \}$ obtained by solving the FOM \eqref{eq:param_sample_pb}, by selecting the first $n$ left singular vectors of the matrix stacking all the snapshots. In particular, prescribed a given tolerance $\varepsilon_{POD}$ on the sum of the discarded modes (quantified by the sum of the squares of the corresponding singular values), a resulting dimension $n=n(\varepsilon_{POD})$ of the POD space is obtained. Note that the RB matrix $V$ does not depend on $\bs{\mu}$ or time. We denote by $\bs{\Xi} = \{ \bs{\mu}^i, i=1,\ldots, n_s \}$ the train sample of size $n_s$ of parameter values upon which the snapshots are computed.

To obtain the ROM that, for any given $\bs{\mu} \in \mathcal{P}$, allows us to compute the RB coordinates, we project the initial datum onto the RB subspace, getting $\bb{y}_n(0) = V^T\bb{y}_0$, and  introduce the residual of the FOM evaluated on the ROM approximation, \vspace{-0.1cm}
\[
\bb{r}({\bb{y}}_n(t, \bs{\mu}); \bs{\mu}) = \bb{F}(\bs{\mu}) -   {M}  V \dot{\bb{y}}_n(t, \bs{\mu}) -  {A}(\bs{\mu}) V \bb{y}_n(t, \bs{\mu});  \vspace{-0.1cm}
\]
then, we impose that the residual is orthogonal to the subspace spanned by the $n$ selected basis functions, that is,  \vspace{-0.1cm}
\[
V^T \bb{r}({\bb{y}}_n(t, \bs{\mu}); \bs{\mu}) = V^T \left(
\bb{F}(\bs{\mu}) -   {M}  V \dot{\bb{y}}_n(t, \bs{\mu}) -  {A}(\bs{\mu}) V \bb{y}_n(t, \bs{\mu})
\right) = {\bf 0}, \qquad t \in (0,T). \vspace{-0.1cm}
\]
This yields the following $n$-dimensional ROM, \vspace{-0.1cm}
\begin{equation}
    \begin{cases}
        {M}_n  \dot{\bb{y}}_n +  {A}_n(\bs{\mu})\bb{y}_n = \bb{F}_{n}(\bs{\mu})    \, , \, &t \in (0,T) \\
     \bb{y}_n(0) = {V}^T \bb{y}_0
    \end{cases}
    \vspace{-0.1cm}
    \label{eq:param_sample_pb_ROM}
\end{equation}
where ${M}_n = {V}^T  M {V}$, ${A}_n(\bs{\mu}) = V^T {A}(\bs{\mu}) V$, and $\bb{F}_{n}(\bs{\mu})  = V^T \bb{F}(\bs{\mu})$. Note that the similar strategy is used to obtain a ROM for a steady problem, and ultimately yields a steady ROM, provided the snapshots are computed by solving the (steady) FOM for $n_s$ selected input parameters vectors. Note that all the ROM arrays that do not depend on the parameter vector can be assembled and stored once and for all.

The ROM \eqref{eq:param_sample_pb_ROM} can be efficiently solved relying on the time integrator used for the FOM, provided its arrays can be assembled inexpensively, that is, independently of the FOM dimension $N_y$. This requirement is automatically fulfilled provided   $A$ and ${\bf F}$ are affinely parametrized, i.e, \vspace{-0.1cm}
\[
{A}(\bs{\mu}) = \sum_{q=1}^{Q_A} \Theta_q^A(\bs{\mu}) A^q, \qquad \qquad
{\bf f} (\bs{\mu}) = \sum_{q'=1}^{Q_F}  \Theta_{q'}^F(\bs{\mu}) {\bf F}^{q'} \vspace{-0.1cm}
\]
for a set of $\bs{\mu}$-dependent functions $\Theta_q^A(\bs{\mu})$,  $\Theta_{q'}^F(\bs{\mu})$, and $\bs{\mu}$-independent matrices $A^q$, $q=1,\ldots,Q_A$ and vectors  ${\bf F}^{q'}$, $q'=1,\ldots,Q_F$. Under this assumption, we can express \vspace{-0.1cm}
\[
\begin{array}{c}
{A}_n(\bs{\mu}) = V^T {A}(\bs{\mu}) V = \sum_{q=1}^{Q_A} \Theta_q^A(\bs{\mu}) V^T A^q V =  \sum_{q=1}^{Q_A} \Theta_q^A(\bs{\mu}) A_n^q, \\
{\bf F}_n(\bs{\mu}) = V^T {\bf F}(\bs{\mu}) V = \sum_{q'=1}^{Q_F} \Theta_{q'}^F(\bs{\mu}) V^T {\bf F}^{q'} =  \sum_{q'=1}^{Q_F} \Theta_{q'}^F(\bs{\mu}) {\bf F}_n^{q'}, 
\end{array} \vspace{-0.1cm}
\]
so that we can precompute and store all the (low-dimensional) ROM arrays $A^q_n$, $q=1,\ldots,Q_A$ and ${\bf F}^{q'}_n$, $q'=1,\ldots,Q_F$, and retrieve the parameter-dependent ROM operators inexpensively. In the case the affine parameteric assumption is not automatically fulfilled by the original FOM, suitable hyper-reduction strategies must be considered, however potentially limiting the overall computational speedup offered by the ROM. See, e.g., \cite{red_book,Negri2015} for further details.

\subsection{Reduced order modeling of parametrized OCPs}

The reduction strategy we employ to speedup the solution of the OCPs we focus on is based on a projection of the reference, the state, the adjoint and the control variables onto lower dimensional subspaces, following the strategy described for \emph{parametrized} OCPs in \cite[Chapter 12]{red_book}.  The mathematical theory behind the construction of ROMs for linear quadratic OCPs is by now well-developed and has proven to be powerful in many OCPs for parametrized PDEs, among many others we highlight the model-order reduction of fluid flow, advection-diffusion and acoustic problems (see, e.g.,\cite{Manzoni2017,NegriRozzaManzoniQuarteroni2013,NEGRI2015319} ), as well as, more recently, transient problems \cite{Strazzullo2020}. 

\subsubsection{Steady state OCP}

Following the strategy proposed in \cite{NegriRozzaManzoniQuarteroni2013,NEGRI2015319} for the case of steady parametrized OCPs, we build a ROM for the parametrized optimality system \eqref{big_ss} relying on a POD-Galerkin RB method. In particular, we look at  \eqref{big_ss} as a unique parametrized linear system, of the form  \vspace{-0.1cm}
\begin{equation} \label{big_ss_param}
\mathcal{K}(\bs{\mu}) \mathbf{y}(\bs{\mu}) 
= \mathcal{F}(\bs{\mu}), \qquad \qquad 
\mathbf{y}(\bs{\mu})  =  \begin{bmatrix}
\bb{z}^{ss}(\bs{\mu}) \\
\bb{q}^{ss}(\bs{\mu}) \\
\bb{p}^{ss}(\bs{\mu}) \\
\bb{u}^{ss}(\bs{\mu}) 
\end{bmatrix}
\end{equation}
in which the matrix $A$ and the vectors ${\bf F}$ and ${\bf F}_o$ -- appearing as blocks of $\mathcal{K}$ and $\mathcal{F}$, respectively -- depend on the parameter vector $\bs{\mu}$. Here $\mathcal{K}(\bs{\mu}) \in \mathbb{R}^{(N_z+2N_q +N_u) \times (N_z+2N_q +N_u)}$ and $\mathcal{F}(\bs{\mu}) \in \mathbb{R}^{N_z+2N_q +N_u}$. A POD-Galerkin ROM is then obtained by projecting problem \eqref{big_ss_param} onto a (product) space spanned by a set of POD modes for each component. To do this, we first  {\em (i)} solve \eqref{big_ss_param} for $n_s$ parameter vectors sampled from $\mathcal{P}$, and then {\em (ii)} perform POD on each component. Collecting the snapshots onto the following matrices   \[
\begin{array}{c}
S_{qp} =  [\bb{q}^{ss}(\bs{\mu}^1) \, | \, \ldots \, | \,  \bb{q}^{ss}(\bs{\mu}^{n_s})| \, \ldots \, | \,  \bb{p}^{ss}(\bs{\mu}^1) \, | \, \ldots \, | \,  \bb{p}^{ss}(\bs{\mu}^{n_s})], \smallskip \\
S_{z} =  [\bb{z}^{ss}(\bs{\mu}^1) \, | \, \ldots \, | \,  \bb{z}^{ss}(\bs{\mu}^{n_s})], \qquad S_u =  [\bb{u}^{ss}(\bs{\mu}^1) \, | \, \ldots \, | \,  \bb{u}^{ss}(\bs{\mu}^{n_s})]
\end{array}   
\]
the selected POD modes are then used to form the matrices $V_z$, $V_s$ and $V_{qp}$; this latter is a basis matrix used to express both state and adjoint variables. Indeed, basis vectors for the state and the adjoint variables are combined together to form a unique space to approximate both components to ensure the well-posedness of the ROM problem \vspace{-0.1cm}
\begin{equation} 
\underbrace{\mathcal{V}^T\mathcal{K}(\bs{\mu})  \mathcal{V}}_{\mathcal{K}_n(\bs{\mu})} 
 \mathbf{y}_n(\bs{\mu}) 
= \underbrace{\mathcal{V}^T \mathcal{F}(\bs{\mu})}_{\mathcal{F}_n(\bs{\mu})}, \qquad 
\mathbf{y}_n(\bs{\mu}) =\begin{bmatrix}
\bb{z}_n^{ss}(\bs{\mu}) \\
\bb{q}_n^{ss}(\bs{\mu}) \\
\bb{p}_n^{ss}(\bs{\mu}) \\
\bb{u}_n^{ss}(\bs{\mu}) 
\end{bmatrix}, \qquad 
\mathcal{V} = 
\begin{bmatrix}
V_z& \phantom{} & \phantom{} & \phantom{} \\
\phantom{} & V_{qp} & \phantom{} & \phantom{} \\
\phantom{}    &  \phantom{}      & V_{qp}& \phantom{} \\
\phantom{} & \phantom{} & \phantom{}   & V_u
\end{bmatrix};
\end{equation}
a detailed explanation of the construction of the RB spaces in this case can be found, e.g., in \cite{NegriRozzaManzoniQuarteroni2013,NEGRI2015319}. Note that the matrices $V_{qp}$, $V_z$ and $V_u$ have in general a different number of columns, that is, RB spaces of different dimensions are usually generated to approximate the four fields.

\subsubsection{Transient OCP}

In the case of the transient OCP, we consider the simultaneous reduction of the optimality system and the reference dynamics, thus building a low-dimensional approximation for the solution manifold of the OCP and the reference, too. In principle, one could follow the same reduction steps as for the steady-state problem by introducing an additional time discretization. The resulting dimensions of the KKT matrix $\tilde{\mathcal{K}} \in \mathbb{R}^{N(N_z+2N_q +N_u) \times N(N_z+2N_q +N_u)} $ quickly become computationally intractable for relatively fine space and time discretization. Furthermore, we face an additional numerically heavy task in performing the POD reduction of the snapshot matrices $\bb{Z},\bb{Q},\bb{P} $ and $\bb{U} $ which are $N$ times wider than the corresponding steady-state ones. As a consequence, we resort to the iterative Quasi-Newton Algorithm (\ref{alg_newton}) to generate snapshot matrices for a set of scenario parameters sampled from $\mathcal{P}$ using the Latin Hypercube Sampling (LHS) technique, and solving  problem \eqref{eq:param_opt_discrete} for $n_s$ parameter vectors sampled from $\mathcal{P}$.

 Then, the generation of the $V_z$,$V_{pq}$ and $V_u$ is performed sequentially in the parameter set without explicitly performing the reduction of the complete snapshot matrices. In order to do this, for each element of the parameter set, we extract the corresponding snapshot matrix and we build the associated reduced basis sequentially. For example, regarding $\bb{Z}$ we extract $n_s$ matrices $\bb{Z}_i$ with dimension $N_z \times N$ and perform the POD reduction at each step updating the associated reduced basis $V_z$. The conceptual steps of the transient reduction are reported in Algorithm (\ref{alg_pod_transient}). The resulting matrices are then used to project the optimality system \eqref{opt_discrete} onto a lower dimensional subspace, with the same ideas of the previous sections. In particular, the same basis matrix $V_{qp}$, built through POD from snapshots of both the state and the adjoint problems, is used to express the ROM approximation of both these fields.  Also in this case the matrices $V_{qp}$, $V_z$ and $V_u$ have in general a different number of columns -- for ease of notation, we do not distinguish among those dimensions when inserting a  subscript $n$ on the RB approximations of the four fields.
 
 The reduced-order approximation of the KKT system  \eqref{eq:param_opt_discrete}   thus reads as follows: 

\begin{equation}
    \begin{cases}
        V_{r}^{\top}MV_{r} \dot{\bb{z}}_n + V_{r}^{\top}A(\bs{\mu})V_{r}\bb{z}_n = V_{r}^{\top}\bb{F}(\bs{\mu})  \, , \, & \quad t \in (0,T)  \\
        \bb{z}_n(0)                = \bb{0} & \smallskip \\ 
        V_{pq}^{\top}\tilde{M} V_{pq} \dot{\bb{q}}_n + V_{pq}^{\top}\tilde{A}(\bs{\mu})V_{pq}\bb{q}_n = V_{pq}^{\top}E\bb{F}(\bs{\mu}) + V_{pq}^{\top}\bb{F}_{o}(\bs{\mu}) + V_{pq}^{\top}BV_{u}\bb{u}_n   \, , \, & \quad t \in (0,T) \\
     \bb{q}_n(0) = \bb{0} & \smallskip \\
    -V_{pq}^{\top}\tilde{M}V_{pq}   \dot{\bb{p}}_n + V_{pq}^{\top}\tilde{A}(\bs{\mu})V_{pq}\bb{p}_n = V_{pq}^{\top}M_{obs}(V_{pq}\bb{q}_n-V_{r}E\bb{z}_n)   \, , \, & \quad t \in (0,T) \\ 
    \bb{p}_n(T)  = \bb{p}^{ss}_{n} & \smallskip  \\
    \Big(\beta V_{u}^{\top} M_u V_{u} + \beta_g V_{u}^{\top} A_u V_{u} \Big)\bb{u}_n + V_{u}^{\top} B^{\top} V_{pq} \bb{p}_n          = \bb{0}  \, , \, & \quad t \in (0,T).
    \end{cases}
\label{ocp_rom_transient}
\end{equation}
We highlight that the transient reduction contains a nested steady-state reduction due to the presence of the final adjoint condition which needs the solution of the steady-state problem.  The problems appearing in \eqref{ocp_rom_transient} are sequentially solved, relying on the modified Newton algorithm \eqref{alg_newton}; in this latter, each query to a FOM problem must be replaced by the  query to the corresponding ROM. A similar strategy had been considered in \cite{manzoni2019certified} to address the RB approximation in the simpler case of   steady, linear parametric OCPs, in which parameters were used to define a (low-dimensional) control function. 

\begin{algorithm}
  \begin{algorithmic}[1]
  \State $\text{Snapshot generation from FOM}$
    \For{$\bs{\mu} \,\, \text{in} \,\,  \bs{\Xi}$} \Comment{in parallel}
        \State    $\bb{z},\bb{q},\bb{p},\bb{u} \gets \text{Solve full order OCP with parameter vector } \bs{\mu} \text{ using Algorithm (\ref{alg_newton}) }$ 
        \State $  {Z}, {Q}, {P}, {U}  \gets \bb{z},\bb{q},\bb{p},\bb{u}  $ \Comment{Assemble snapshots matrix}
    \EndFor
    \State $\text{ Reduced basis construction}$
    \For{$i=1:n_s$}
    \State $V_z    = \text{POD}([V_z, {Z}_i],\text{tol})$    
    \State $V_{pq}  = \text{POD}([V_{pq}, {Q}_i, {P}_i],\text{tol})$
    \State $V_u     = \text{POD}([V_u, {U}_i],\text{tol})$
    \EndFor
    
  \end{algorithmic}
  \caption{POD reduction of the parabolic OCP  }
  \label{alg_pod_transient}
\end{algorithm}

\section{Numerical experiments}
\label{NUM}
In this Section, we present several numerical tests both in the transient and steady-state regimes, and for different layouts of the cloak and the obstacle. After numerically checking that solving the optimality system (\ref{opt_discrete}) allows to converge to the steady-state solution at the final instant, we investigate the performance of the reduction technique in capturing the dynamics of the \emph{parametrized} OCP arising from the cloaking objective. 

\begin{figure}[h!]
\centering
\subfigure{\includegraphics[width=0.425\textwidth]{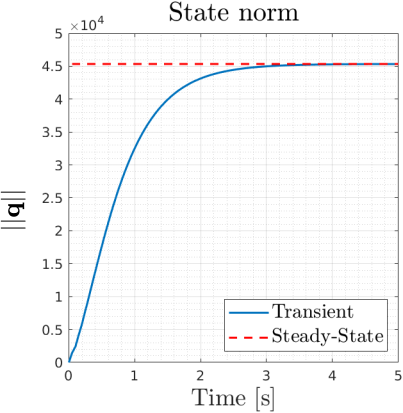}}
\subfigure{\includegraphics[width=0.425\textwidth]{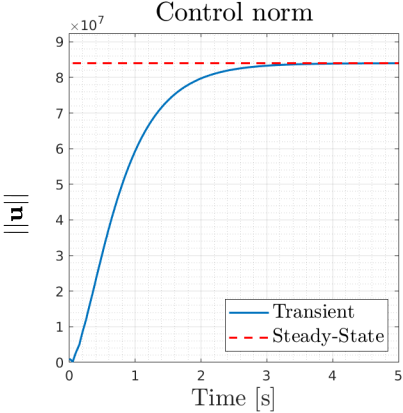}}
\subfigure{\includegraphics[width=0.425\textwidth]{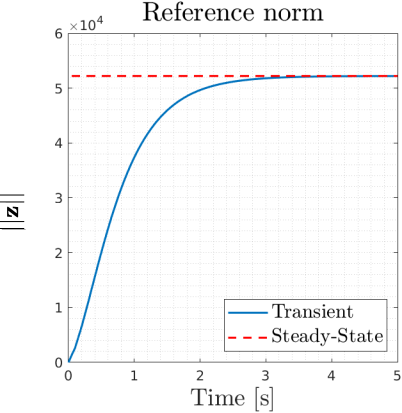}}
\subfigure{\includegraphics[width=0.425\textwidth]{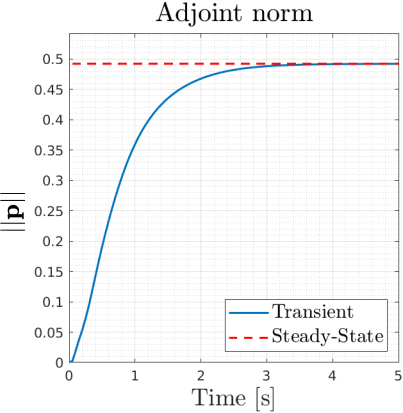}}
\caption{Convergence in the $L^2$-norm of the optimal dynamics of the transient regime to its steady-state counterpart. Note that the final time $T$ is chosen sufficiently large to reach the reference steady-state.}
\label{convergence_fom}
\end{figure}

Before considering ROM performances, we numerically test that replacing the adjoint condition with $p(\x,T)=p^{ss}(\x)$ in the optimality system actually allows the dynamics of the optimal variables to converge to their respective steady-state solution. In order to assess convergence, we use the $L^2$ norm in $\Omega$. Recall that the $L^2$ norm of a function $f \in L^2(\Omega)$ is defined as $\norm{f}_{L^2(\Omega)} = \sqrt{\int_{\Omega} f(\x,t)^2 d\x}$ and is a function of time for time-dependent variables, while it is obviously a constant for steady-state variables. Thus, we want to confirm that 
\begin{equation*}
\begin{array}{c}
    \lim_{T \to \infty} \norm{q(\cdot,T)-q^{ss}}_{L^2(\Omega)} = 0 , \\
     \lim_{T \to \infty} \norm{p(\cdot,T)-p^{ss}}_{L^2(\Omega)} = 0 , \\
     \lim_{T \to \infty} \norm{u(\cdot,T)-u^{ss}}_{L^2(\Omega_c)} = 0, 
     \end{array}
\end{equation*}
that is, for a sufficiently large time horizon $T$, the transient optimal solution converges to the steady-state ones at the final time $T$. Therefore, we solve the parabolic OCP using Algorithm (\ref{alg_newton}) with the parameters set $\bs{\mu} = [\,3.5\,\, 10^4 \,\, 0\,]$  and we check the convergence at their steady-state value. In Figure \ref{convergence_fom} the time history of the $L^2$ norms of the reference and of the optimal variables is shown. The controlled optimal dynamics converge to the steady-state value. Note that we have chosen a sufficiently large time interval $T$ to let the uncontrolled reference reach its steady-state value.

\subsection{Steady-state OCP reduction}

\begin{figure}[h!]
\centering
\subfigure{\includegraphics[width=0.3\textwidth]{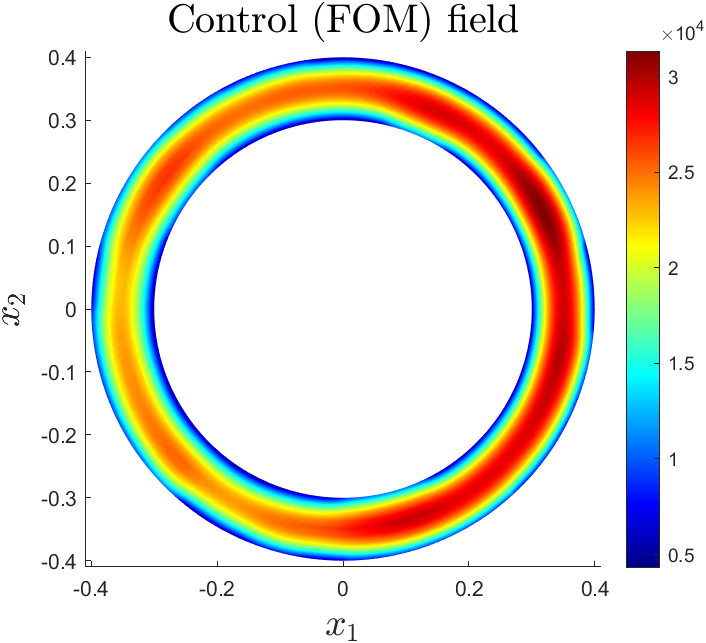}}
\subfigure{\includegraphics[width=0.3\textwidth]{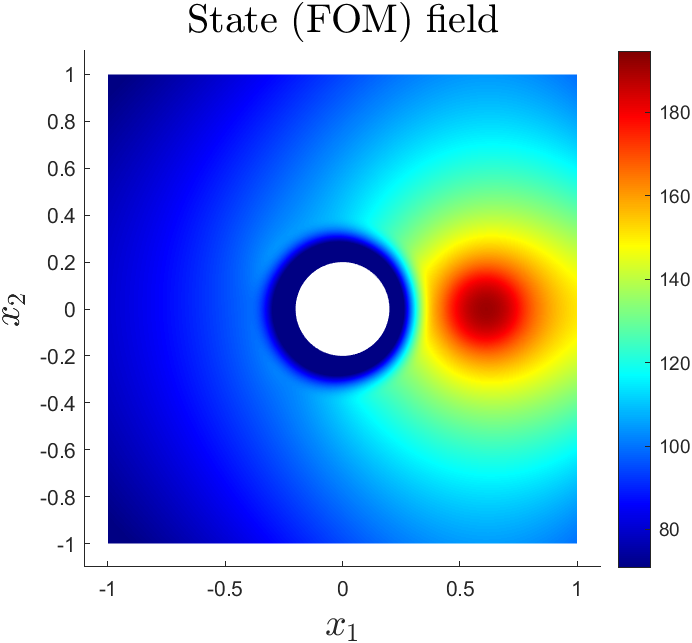}}
\subfigure{\includegraphics[width=0.3\textwidth]{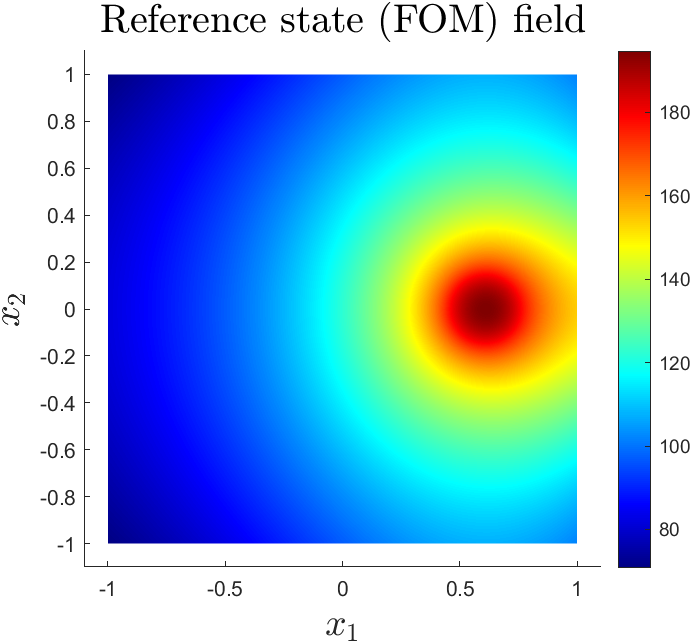}}
\subfigure{\includegraphics[width=0.3\textwidth]{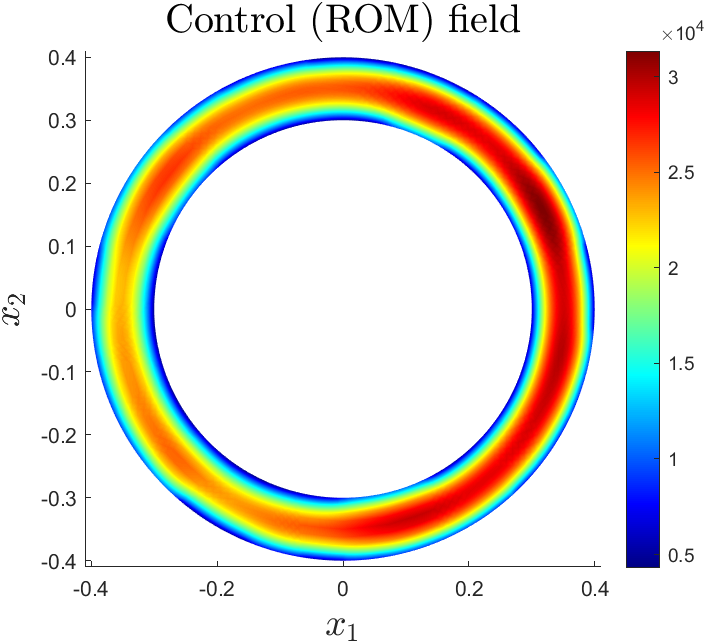}}
\subfigure{\includegraphics[width=0.3\textwidth]{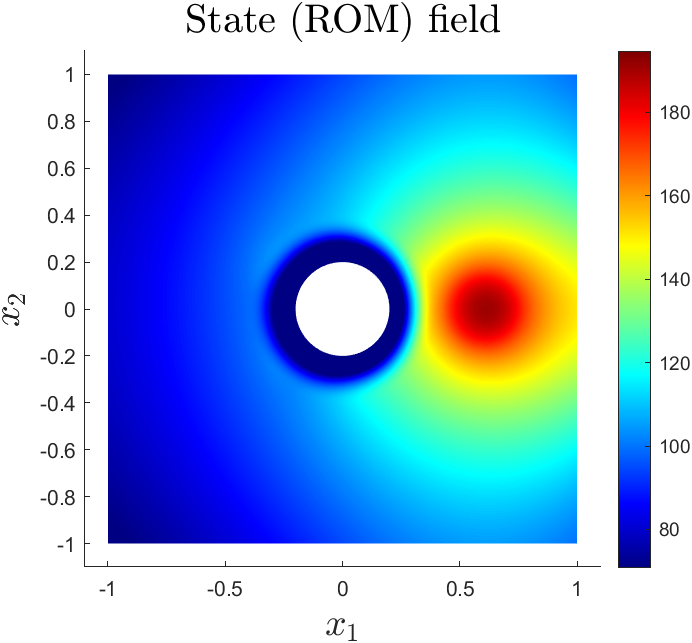}}
\subfigure{\includegraphics[width=0.3\textwidth]{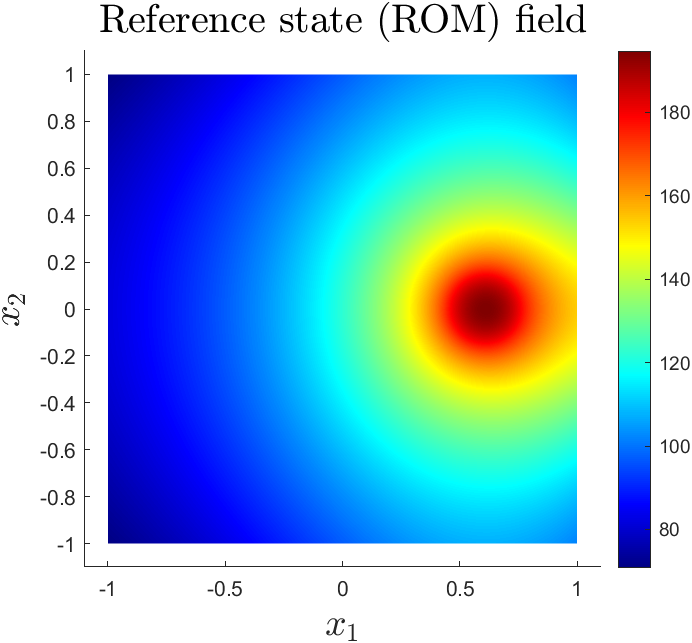}}
\caption{Reduction of the parametrized steady-state OCP, computational results for $\bs{\mu}_{t1} = [\,3.5 \,\, 10^4 \,\, 0 \,] $.}
\label{ss_tc1}
\end{figure}

\begin{figure}[h!]
\centering
\subfigure{\includegraphics[width=0.3\textwidth]{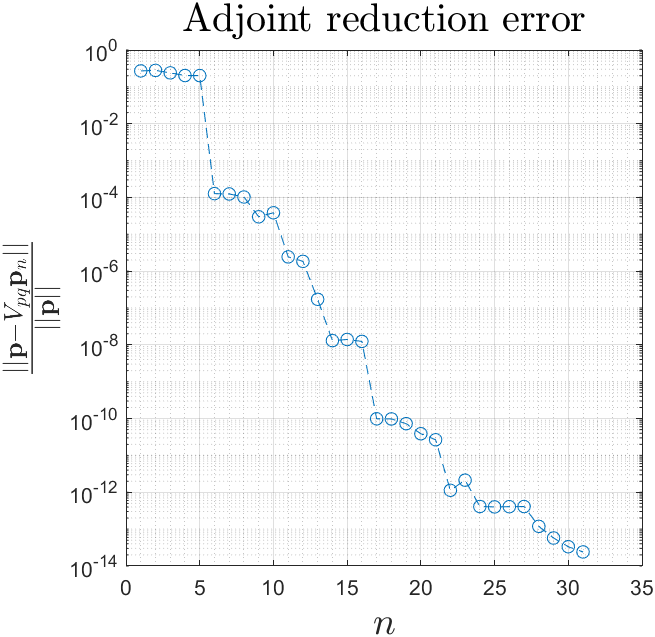}}
\subfigure{\includegraphics[width=0.3\textwidth]{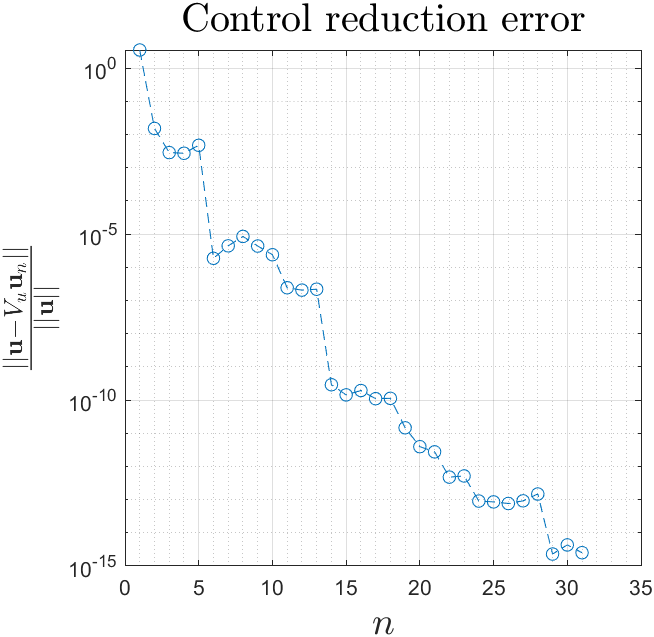}}
\subfigure{\includegraphics[width=0.3\textwidth]{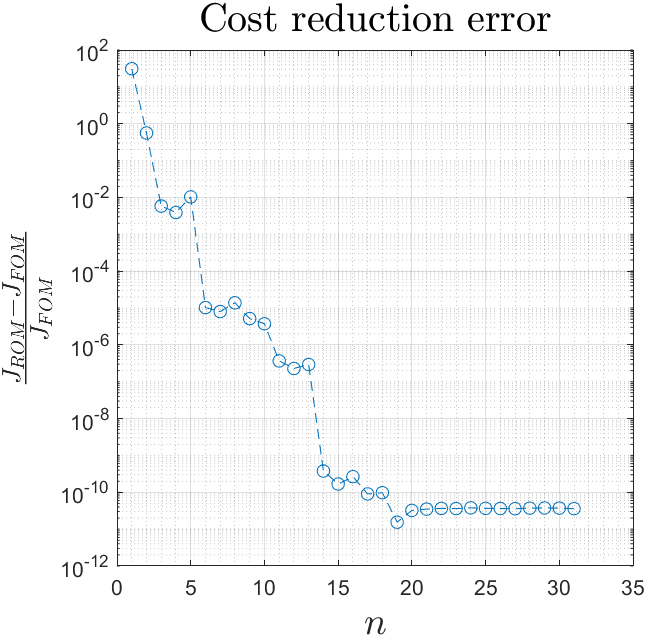}}
\subfigure{\includegraphics[width=0.3\textwidth]{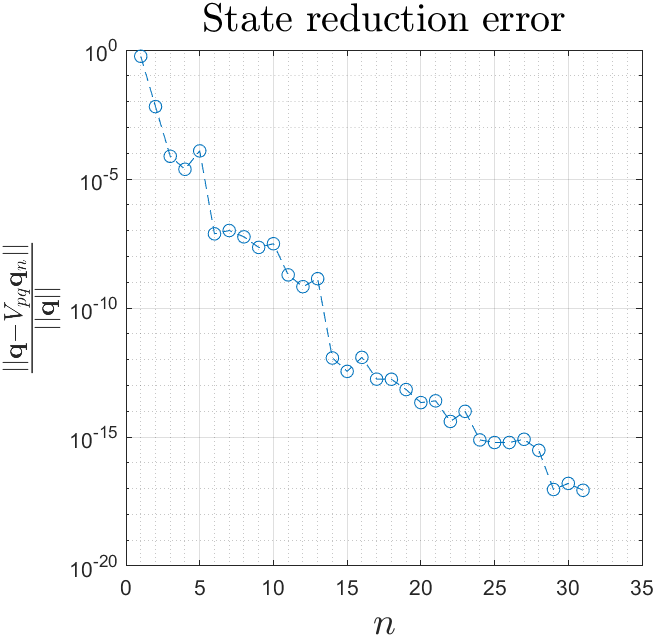}}
\subfigure{\includegraphics[width=0.3\textwidth]{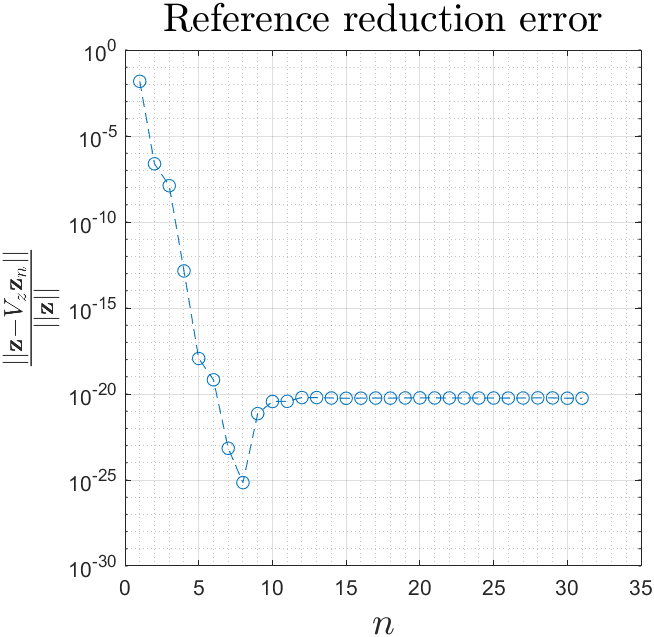}}
\subfigure{\includegraphics[width=0.3\textwidth]{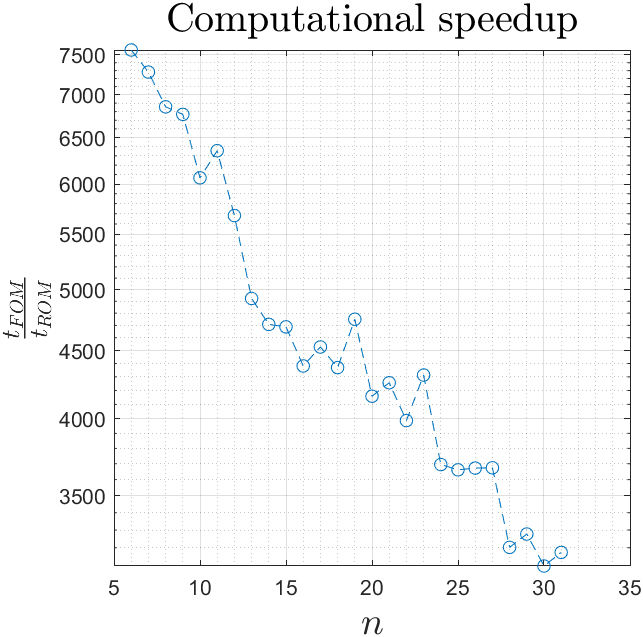}}
\caption{Relative errors between the ROM and FOM solutions. The ROM is able to achieve a computational speedup of $3170$ times without losing in accuracy of the optimal and reference solutions. Solution time for the projected linear system (\ref{big_ss}) is 0.49 milliseconds on a standard laptop computer while it takes 1.57 seconds to solve the FOM.}
\label{conv_ss_tc1}
\end{figure}

We can now proceed to assess the POD-Galerkin ROM for the steady-state OCP. A parameter set of $n_{s}=50$ snapshots is considered using Latin-Hypercube sampling of the parameter space $\mathcal{P} =\left[1,5\right] \times [ \, 5 \cdot 10^2 , 1.5 \cdot 10^4 \,  ] \times \left[ 0 , 200 \right]$. Then, for the sake of visualization, we test the performance of the ROM on $\bs{\mu}_{t1} = [\,3.5 \,\, 10^4 \,\, 0\,] $ and $\bs{\mu}_{t2} = [\,3.5 \,\, 10^4 \,\, 100 \,] $. The dimension of the FOM space discretization is $N_q=17504$, $N_z=18721$, $N_u=3348$ for the state, reference and control variables, respectively. Hence, the sparse linear system (\ref{big_ss}) has dimension $N_z+2N_q+N_u = 57077$. The steady-state problem is solved for each parameter in the set of snapshots using direct methods, 
 while for the assembling of the FEM matrices and the construction of the POD-Galerkin ROM we rely on \texttt{RedbKit} \cite{redbKIT}. Figure \ref{ss_tc1} shows that the FOM and ROM results are practically indistinguishable. Furthermore, the optimal control formulation of the cloaking problem is able to closely follow the reference field. Control parameters $\beta = 10^{-7}$ and $\beta_g = 10^{-8}$ are used in this case. The gradient weighting term $\beta_g$ allows to exploit the whole control domain so that the control intensity is evenly distributed in the control region. When the case $\beta_g=0$ is considered, the control intensity is concentrated in a thin layer close to the boundary of the control region. The numerical results of such a case are not presented for the sake of brevity. The POD algorithm selects $n_q = 31$, $n_z=12$ and $n_u=11$ basis functions with a final computational speedup of roughly $3170$ with respect to the FOM. As a consequence, the ROM is able to compute the optimal solution for a new set of input parameters in 0.49 milliseconds. The relative accuracy of the ROM is shown in Figure \ref{conv_ss_tc1} where, for a sufficiently rich basis, we are able to reach approximation errors close to machine epsilon.

\subsection{Transient OCP reduction}
\begin{figure}[h!]
\centering
\subfigure{\includegraphics[width=0.3\textwidth]{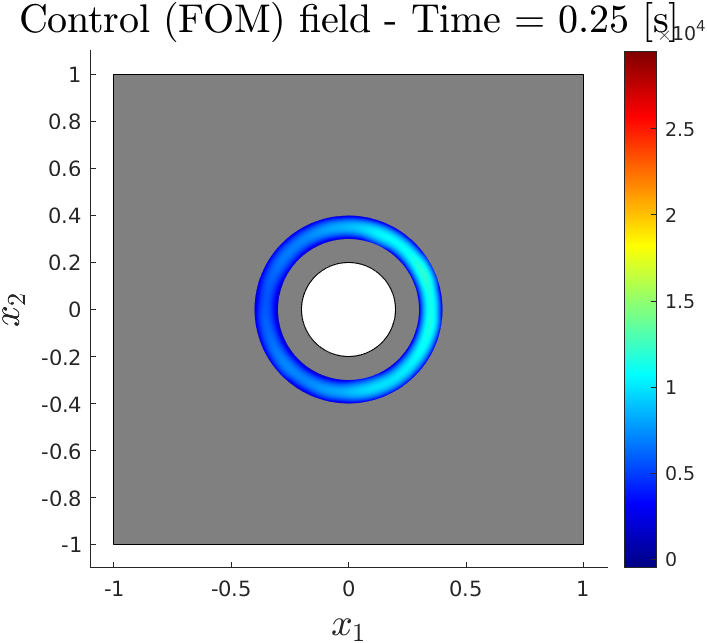}}
\subfigure{\includegraphics[width=0.3\textwidth]{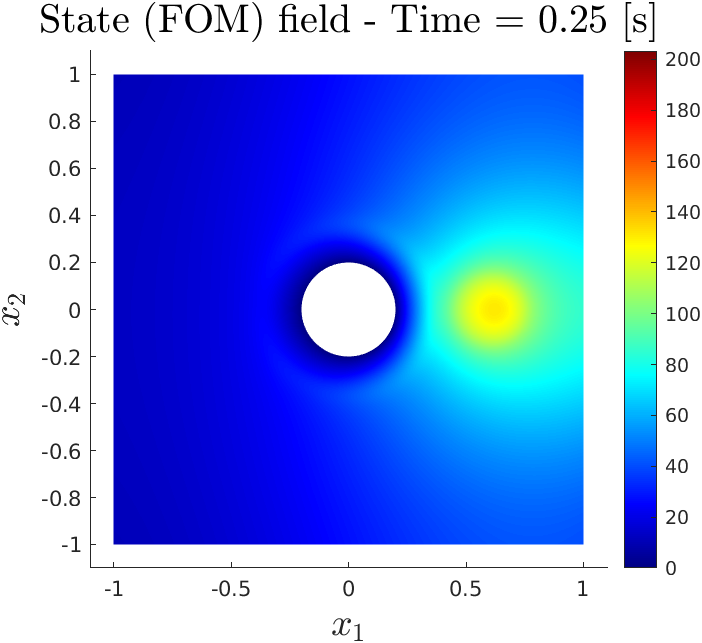}}
\subfigure{\includegraphics[width=0.3\textwidth]{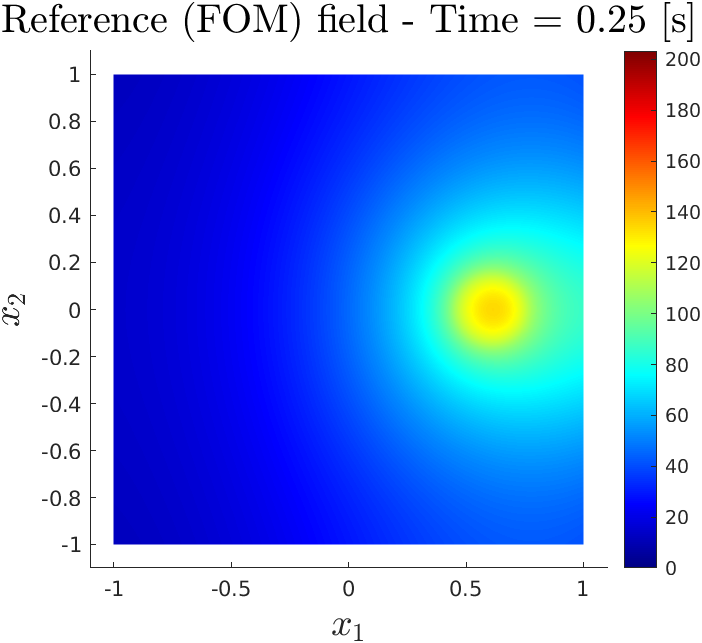}}
\subfigure{\includegraphics[width=0.3\textwidth]{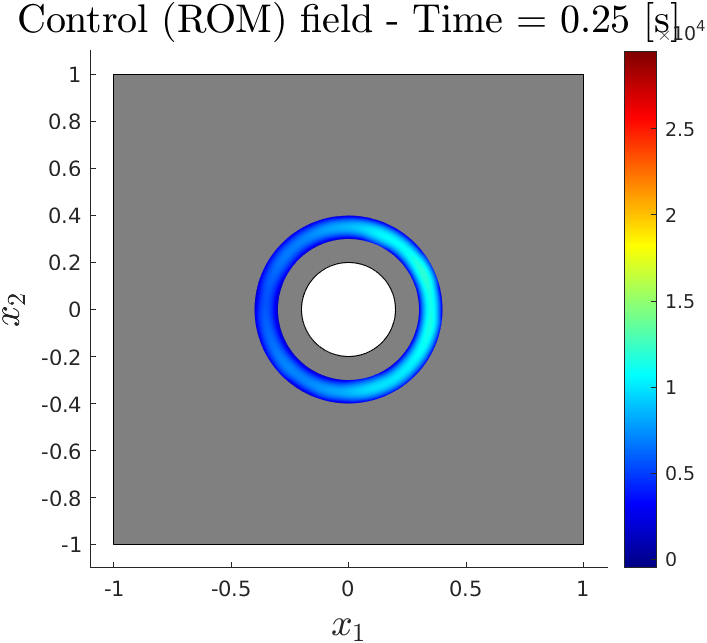}}
\subfigure{\includegraphics[width=0.3\textwidth]{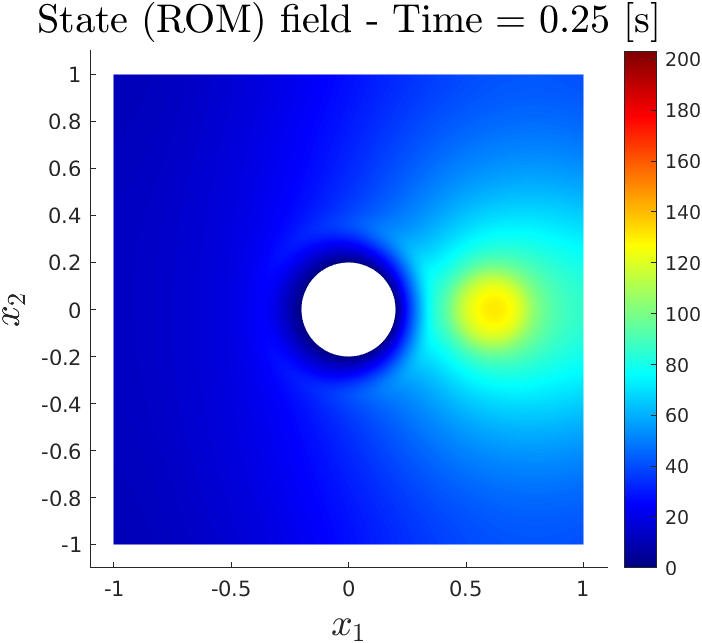}}
\subfigure{\includegraphics[width=0.3\textwidth]{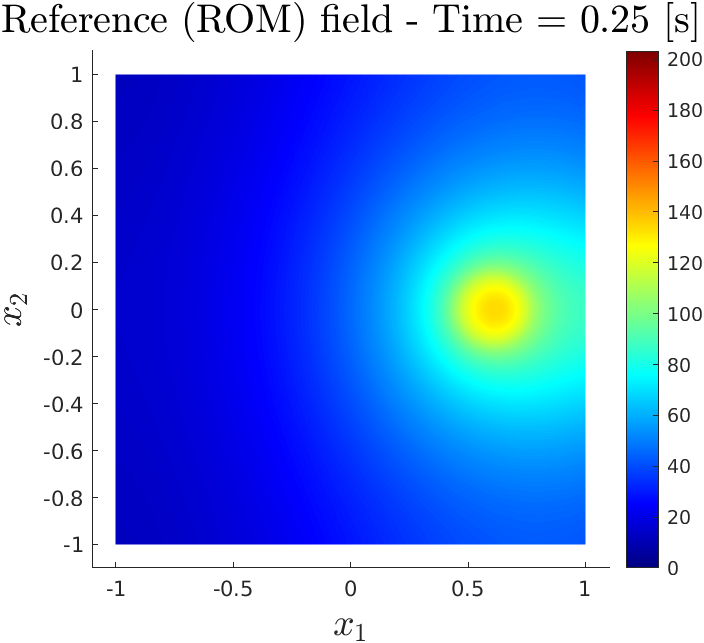}}
\caption{Reduction of the parametrized transient OCP, computational results for $\bs{\mu}_{t1} = [3.5 \, 10^4 \, 0] $ at time instance $t=0.25(\text{s})$.}
\label{t_tc1_frame_1}
\end{figure}

\begin{figure}[h!]
\centering
\subfigure{\includegraphics[width=0.3\textwidth]{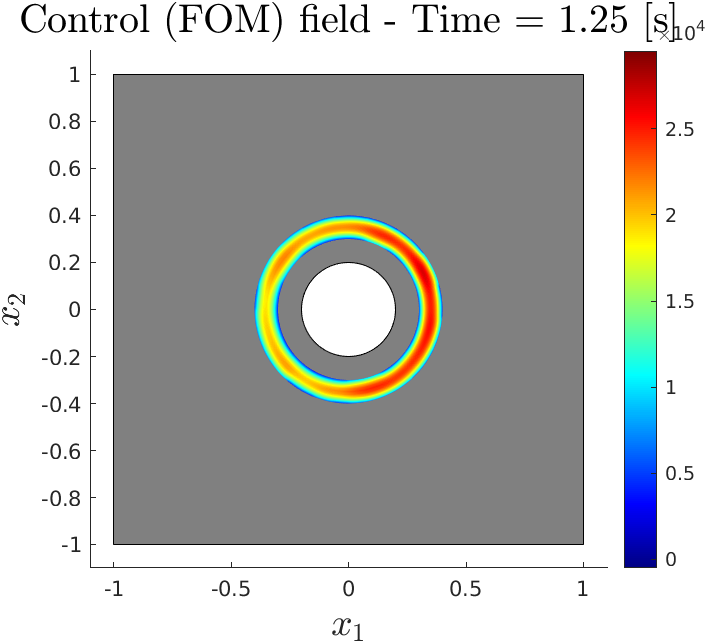}}
\subfigure{\includegraphics[width=0.3\textwidth]{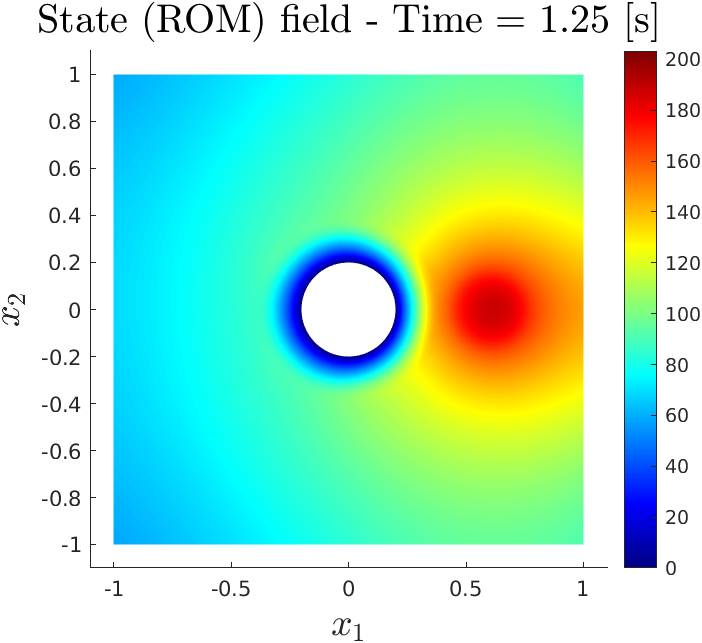}}
\subfigure{\includegraphics[width=0.3\textwidth]{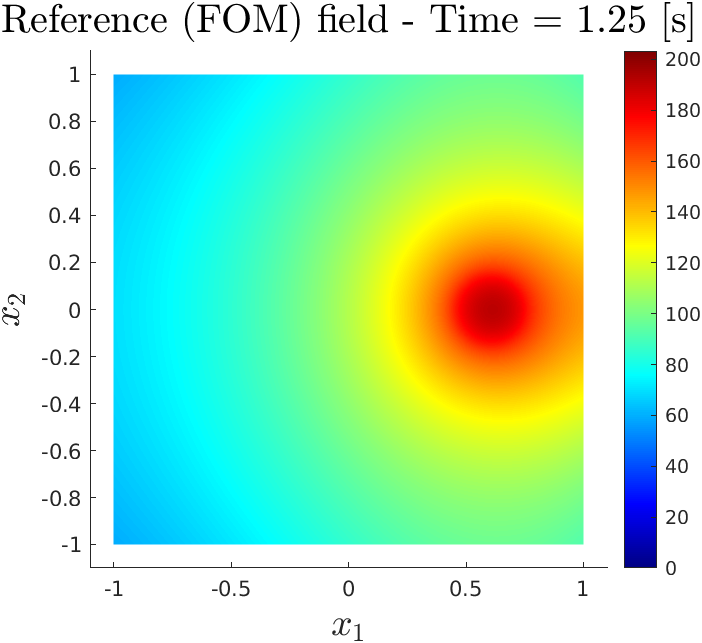}}
\subfigure{\includegraphics[width=0.3\textwidth]{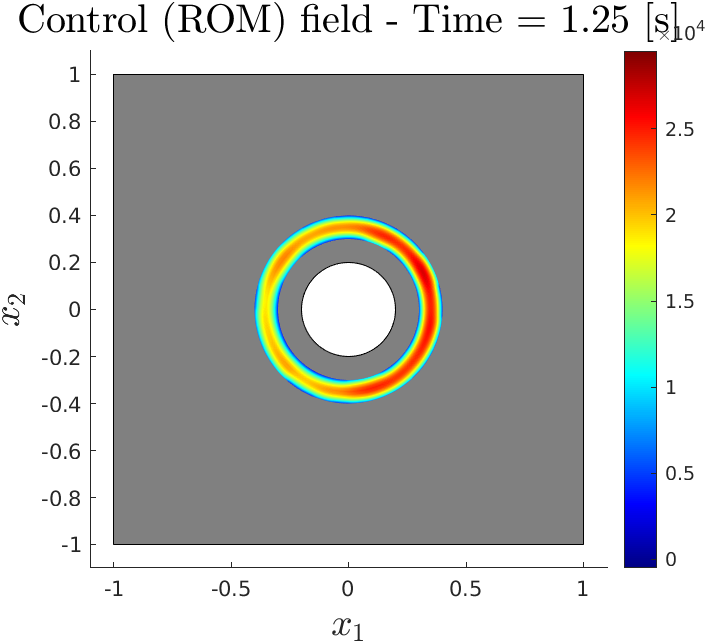}}
\subfigure{\includegraphics[width=0.3\textwidth]{images/T_tc1/state_ROM_field_26.png}}
\subfigure{\includegraphics[width=0.3\textwidth]{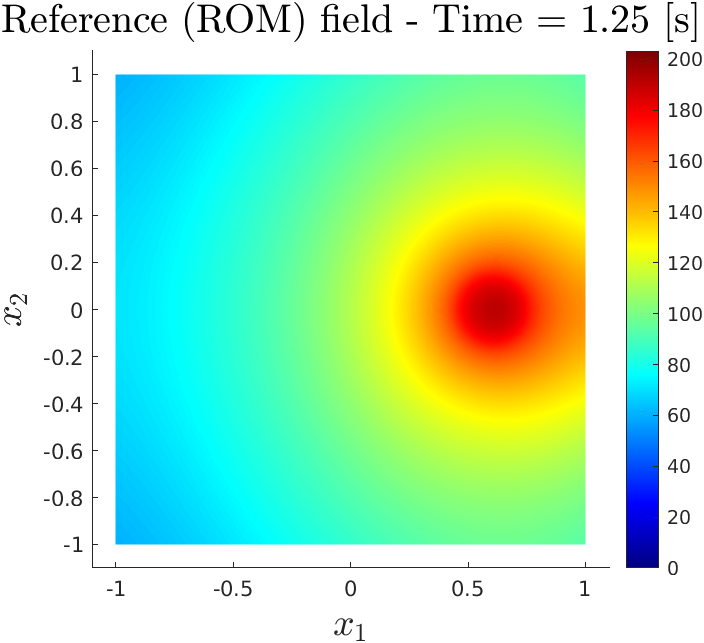}}
\caption{Reduction of the parametrized transient OCP, computational results for $\bs{\mu}_{t1} = [3.5 \, 10^4 \, 0] $ at time instance $t=1.25(\text{s})$.}
\label{t_tc1_frame_2}
\end{figure}

\begin{figure}[h!]
\centering
\subfigure{\includegraphics[width=0.3\textwidth]{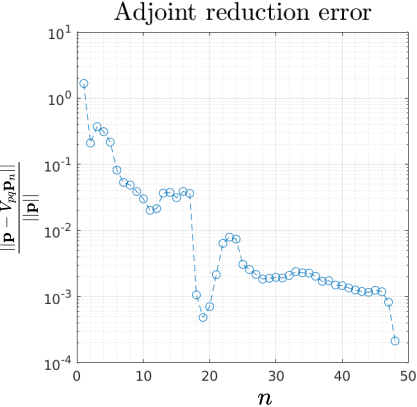}}
\subfigure{\includegraphics[width=0.3\textwidth]{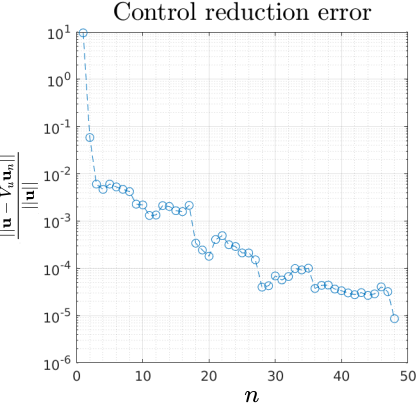}}
\subfigure{\includegraphics[width=0.3\textwidth]{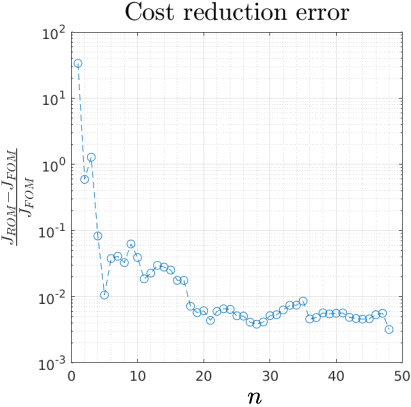}}
\subfigure{\includegraphics[width=0.3\textwidth]{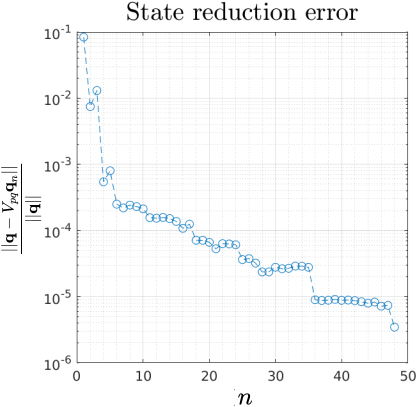}}
\subfigure{\includegraphics[width=0.3\textwidth]{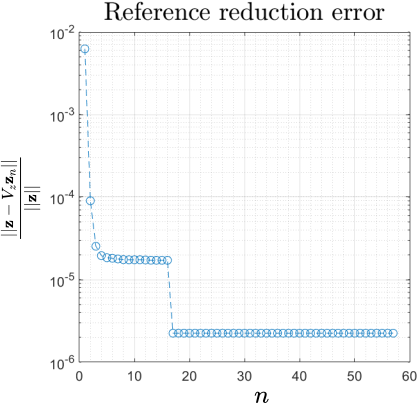}}
\subfigure{\includegraphics[width=0.3\textwidth]{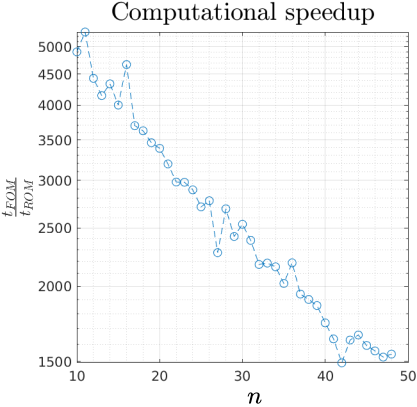}}
\caption{Relative errors between the ROM and FOM solutions for the transient problem. The ROM is able to achieve a computational speedup of $1550$ times without losing in accuracy of the optimal and reference solutions.}
\label{conv_t_tc1}
\end{figure}

We are now ready to tackle the reduction of the transient cloaking problem, that is we use Algorithm (\ref{alg_newton}) to construct a POD-Galerkin ROM for the parametrized parabolic OCP (\ref{ocp_formulation}) with the additional guarantee that the optimal variables of the transient problem converge to their respective steady-state solution as numerically proved in Figure \ref{convergence_fom}. Due to the huge memory needed by a fine discretization in space and time of the transient problem, we limit the snapshot dataset to $n_s=25$ keeping the same parameters range as for the steady-state problem. The discrete (in time and space) version of the optimality system (\ref{opt_discrete}) is still a linear system and could be theoretically solved in one-shot \cite{herzog,Strazzullo2020}, that is the reference and optimal variables in space and time could be computed as the solution of a single yet huge linear system. Discretizing the time interval in $N=100$ time steps, we would end up with a linear system with roughly $N(N_z+2N_q+N_u)\approx 5.5\cdot10^6$ variables, that is intractable for the computational resources of a standard laptop computer. Instead, we solve both the FOM and the ROM using the iterative Algorithm \eqref{alg_newton} and the sequential reduction steps detailed in Algorithm \eqref{alg_pod_transient}. The reduction algorithm selects $n_{q}=48$, $n_{z}=10$ and $n_{u}=20$ reduced basis with a final computational speedup of roughly $1550$ with respect to the FOM solution. The solution of the FOM takes 186 seconds while the ROM needs only 0.11 seconds. A sequence of snapshots taken at time instances $t=0.25 (\text{s})$ and $t=1.25 (\text{s})$ is shown in Figures \ref{t_tc1_frame_1}, \ref{t_tc1_frame_2}. It is clear that the FOM achieves thermal cloaking of the obstacle in the transient domain and that the reduced basis are rich enough to closely follow the evolution both of the reference and of the optimal dynamics. The accuracy of the ROM is confirmed by the convergence results in Figure \ref{conv_t_tc1}. For the transient problem we use the norm in $L^2(0,T,L^2(\Omega))$ to measure errors between the FOM and ROM solutions. Recall that for a function $f \in L^2(0,T,L^2(\Omega))$ this norm is defined as :
\begin{equation*}
\norm{f}_{L^2(0,T,L^2(\Omega))}^2 = \int_0^T \norm{f(\cdot,t)}_{L^2(\Omega)}^2 \, dt.
\end{equation*}

\begin{figure}[h!]
\centering
\subfigure{\includegraphics[width=0.3\textwidth]{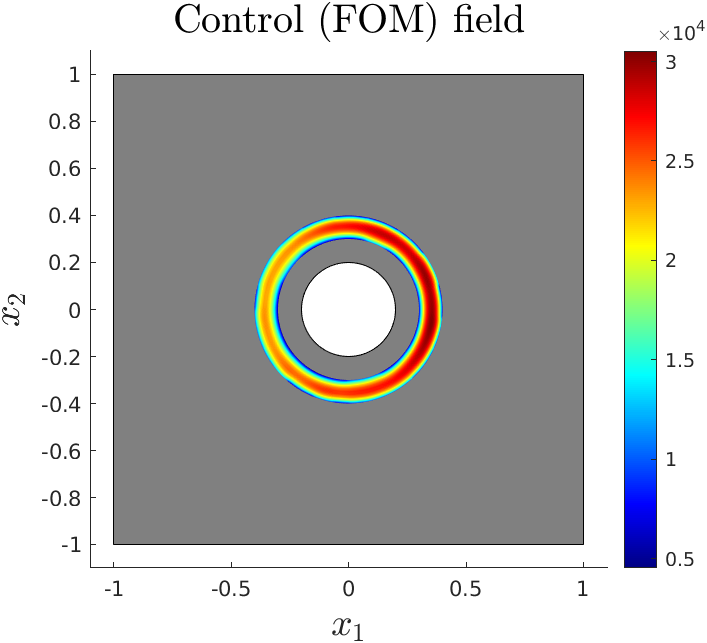}}
\subfigure{\includegraphics[width=0.3\textwidth]{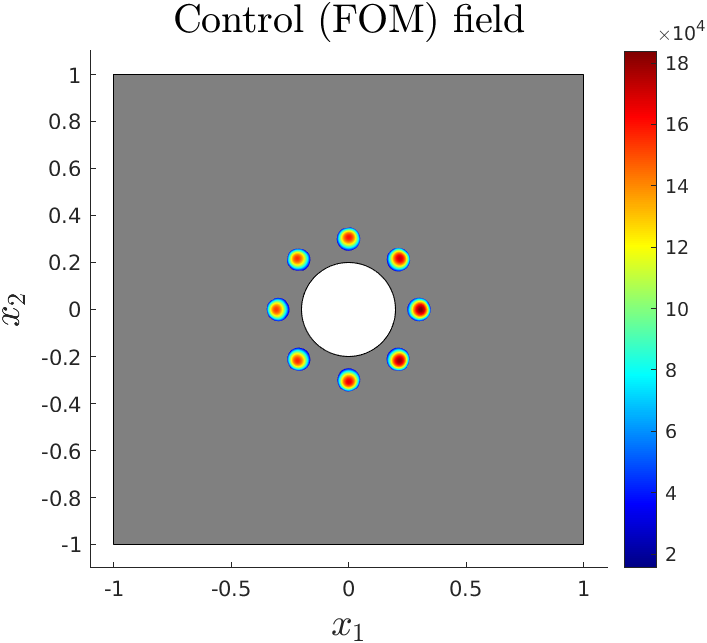}}
\subfigure{\includegraphics[width=0.3\textwidth]{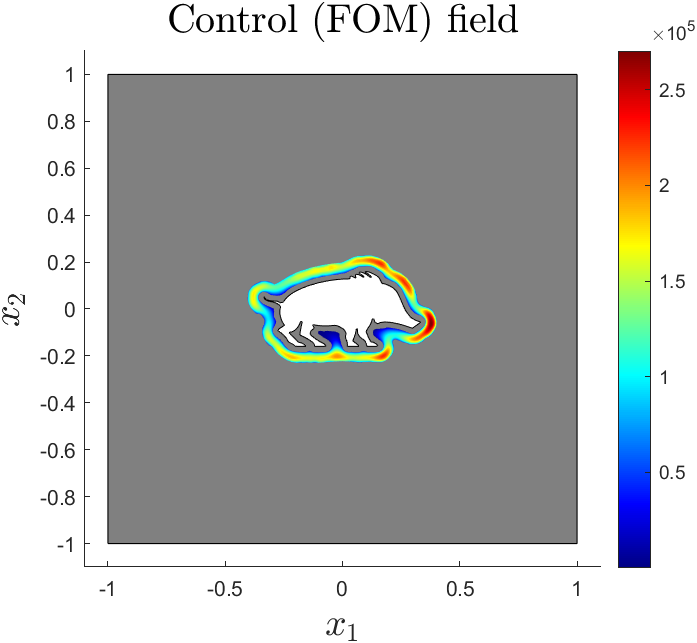}}
\caption{The three cloaking layouts considered in this paper with increasing complexity from left to right. The quantitative assessment of complexity can be seen from the intensity of the control field and from the entity of the tracking errors in Figure \ref{tracking_errors_ss}.}
\label{circles_layout}
\end{figure}

We now turn to the case of a disconnected control domain. In this case, the active sources composing the cloak are disconnected circular domains that surround the obstacle. This layout is inspired by common practical layouts with thermoelectric devices, where it is extremely useful to allow for such flexibility in the control domain. As a test case of cloaking device with detached sources, we consider a set of eight equally spaced heat sources. The control layouts considered in this paper are shown in Figure \ref{circles_layout}, where the second one consists of the disconnected cloak domain while the third one illustrates the complex shape to be treated later. 

\begin{figure}[h!]
\centering
\subfigure{\includegraphics[width=0.3\textwidth]{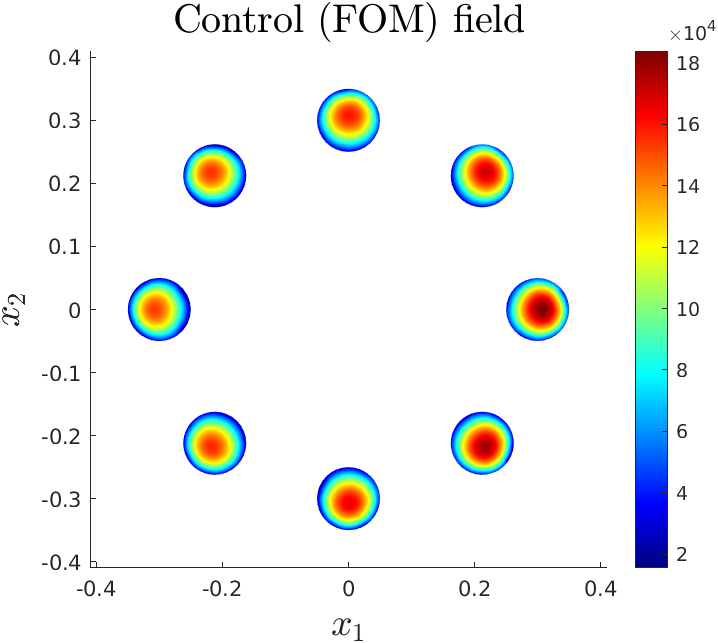}}
\subfigure{\includegraphics[width=0.3\textwidth]{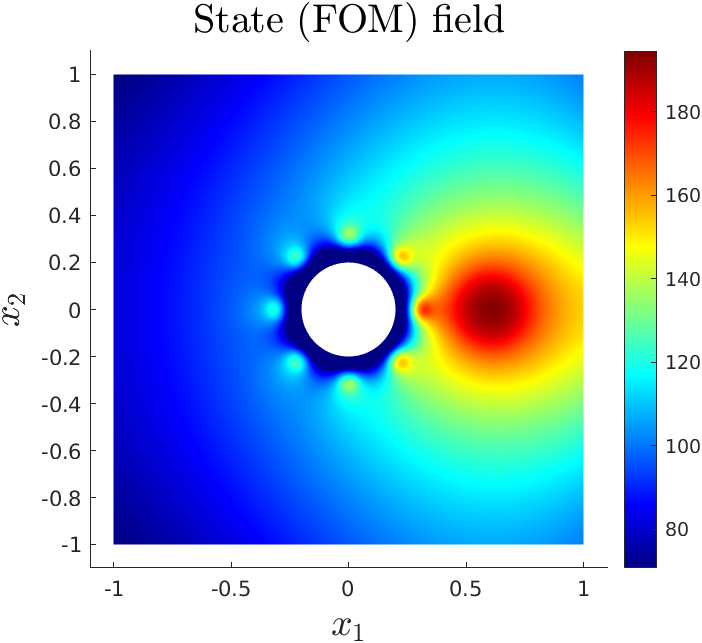}}
\subfigure{\includegraphics[width=0.3\textwidth]{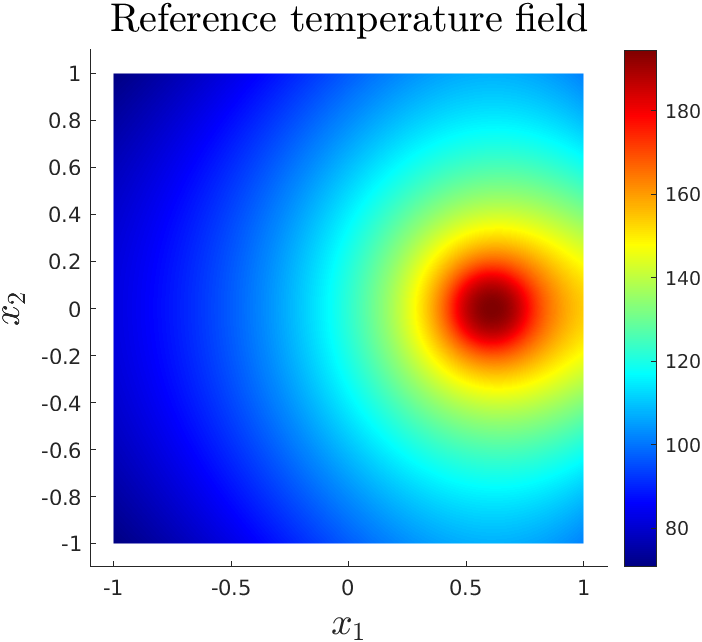}}
\subfigure{\includegraphics[width=0.3\textwidth]{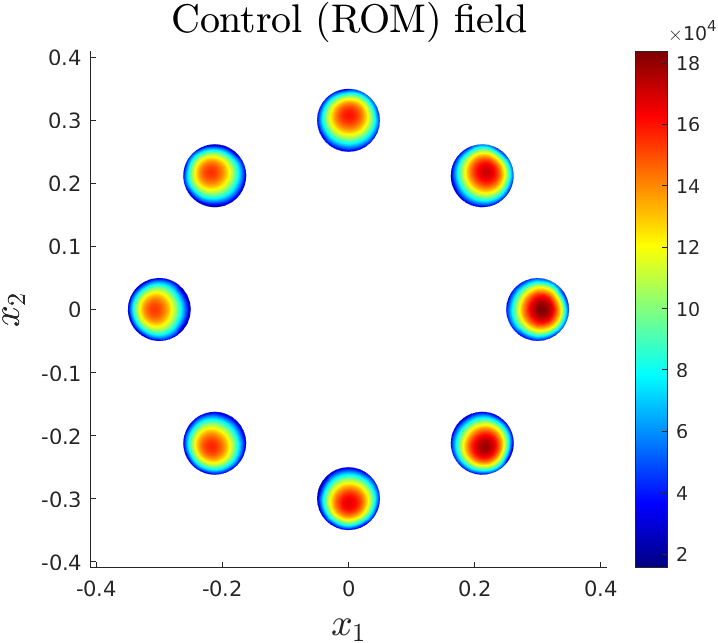}}
\subfigure{\includegraphics[width=0.3\textwidth]{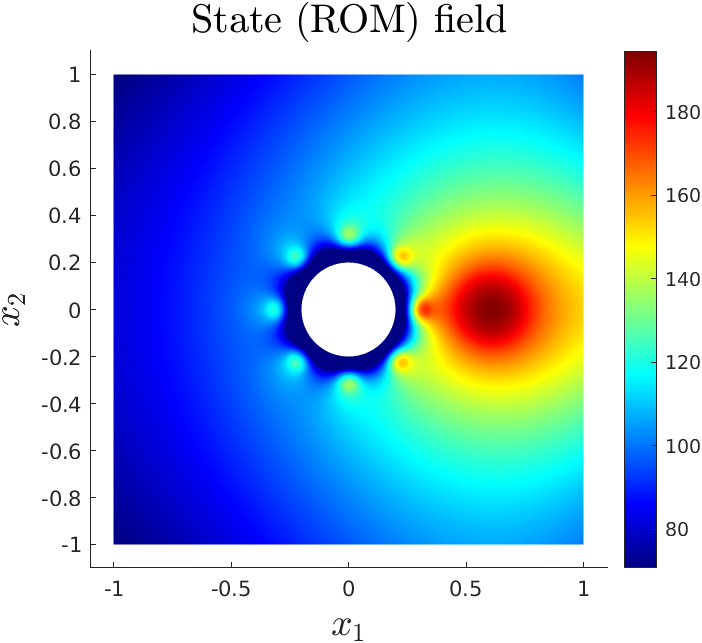}}
\subfigure{\includegraphics[width=0.3\textwidth]{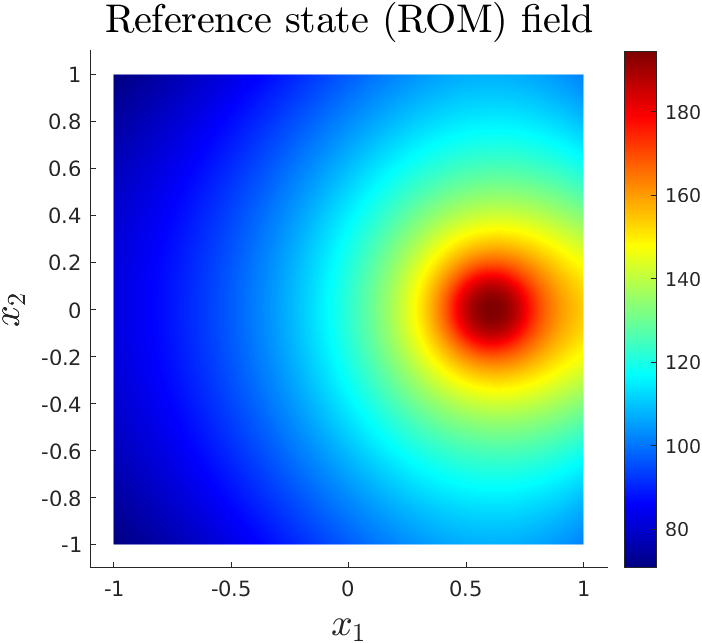}}
\caption{Reduction of the parametrized steady-state OCP with disconnected control domain, computational results for $\bs{\mu}_{t1} = [\,3.5 \,\, 10^4 \,\, 0 \,] $.}
\label{circles_ss}
\end{figure}

\begin{figure}[h!]
\centering
\subfigure{\includegraphics[width=0.3\textwidth]{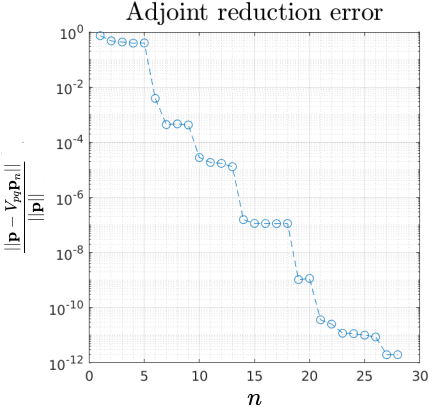}}
\subfigure{\includegraphics[width=0.3\textwidth]{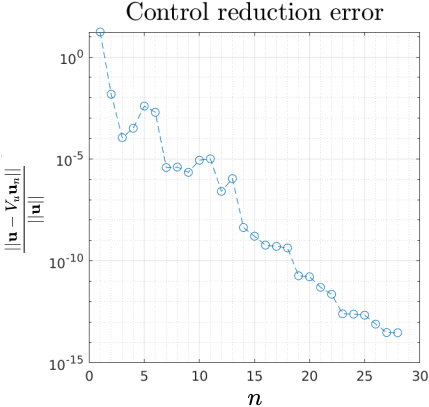}}
\subfigure{\includegraphics[width=0.3\textwidth]{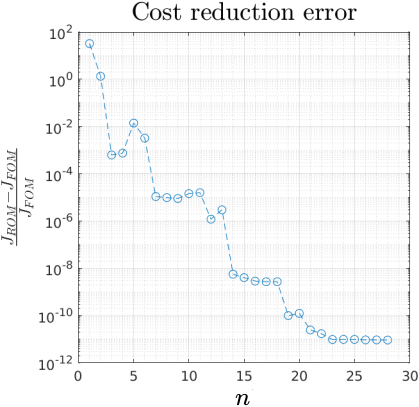}}
\subfigure{\includegraphics[width=0.3\textwidth]{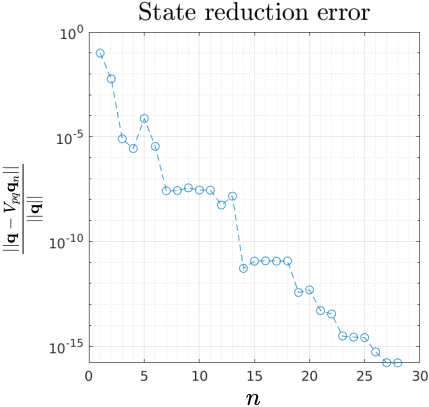}}
\subfigure{\includegraphics[width=0.3\textwidth]{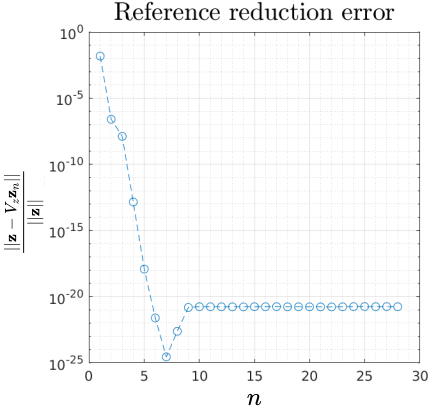}}
\subfigure{\includegraphics[width=0.3\textwidth]{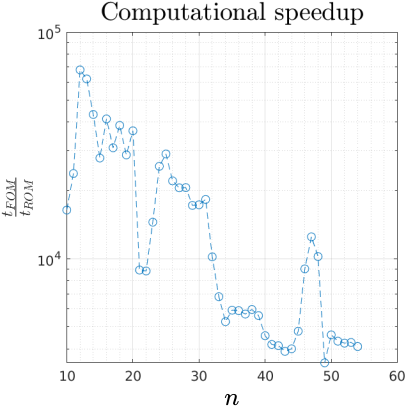}}
\caption{Relative errors between the ROM and FOM solutions for the disconnected control domain. The ROM is able to achieve a computational speedup of $1600$ times without losing in accuracy of the optimal and reference solutions. Solution time for the projected linear system (\ref{big_ss}) is 0.29 milliseconds on a standard laptop computer while it takes 0.46 seconds to solve the FOM.}
\label{circles_ss_ROM}
\end{figure}

\begin{figure}[h!]
\centering
\subfigure{\includegraphics[width=0.3\textwidth]{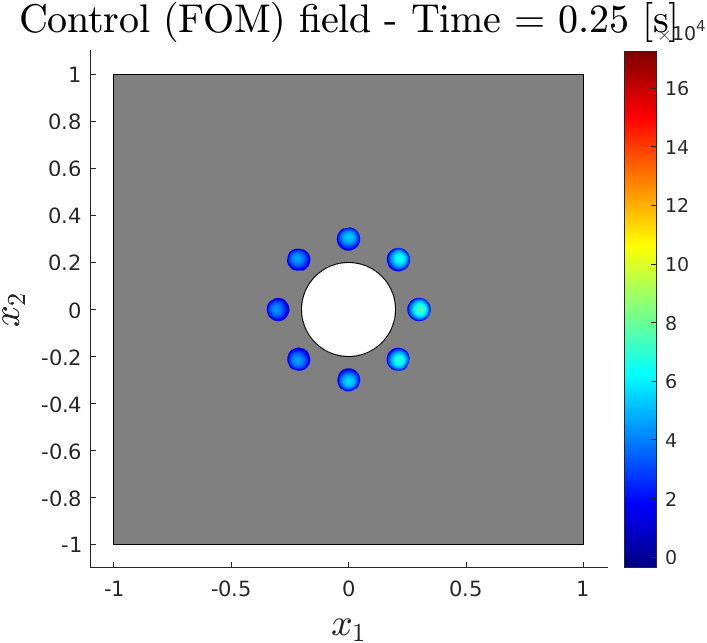}}
\subfigure{\includegraphics[width=0.3\textwidth]{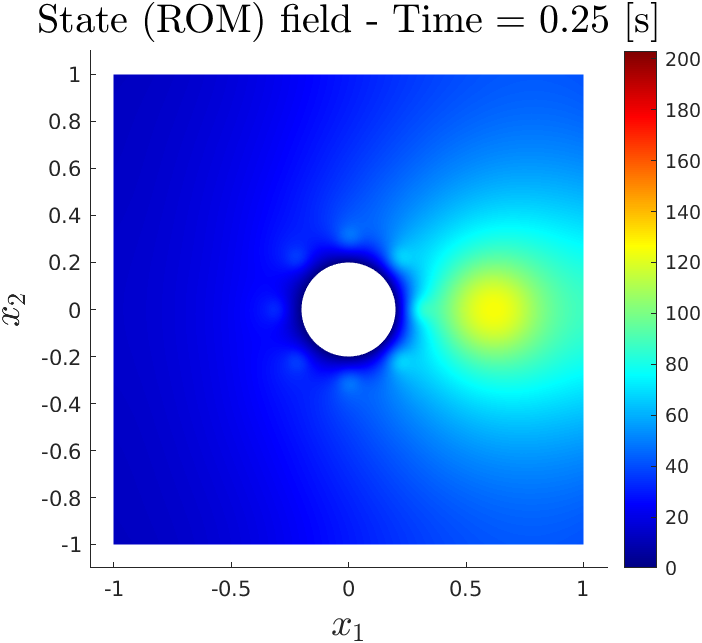}}
\subfigure{\includegraphics[width=0.3\textwidth]{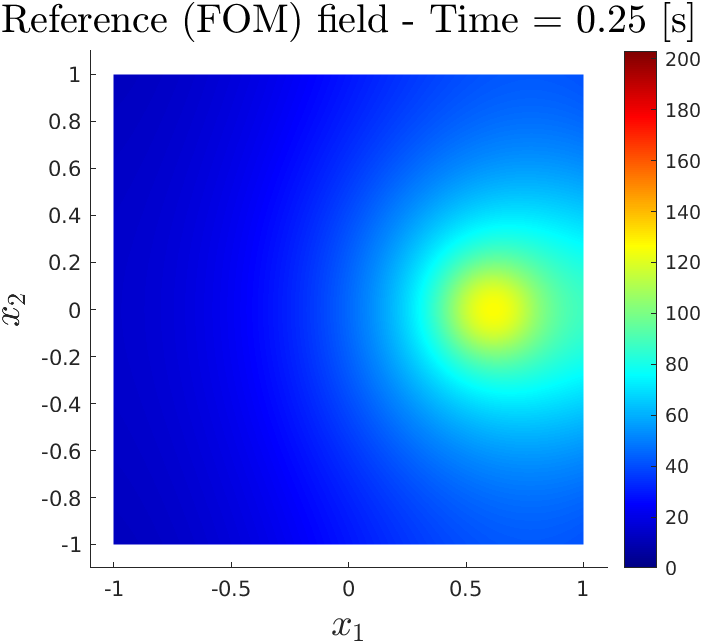}}
\subfigure{\includegraphics[width=0.3\textwidth]{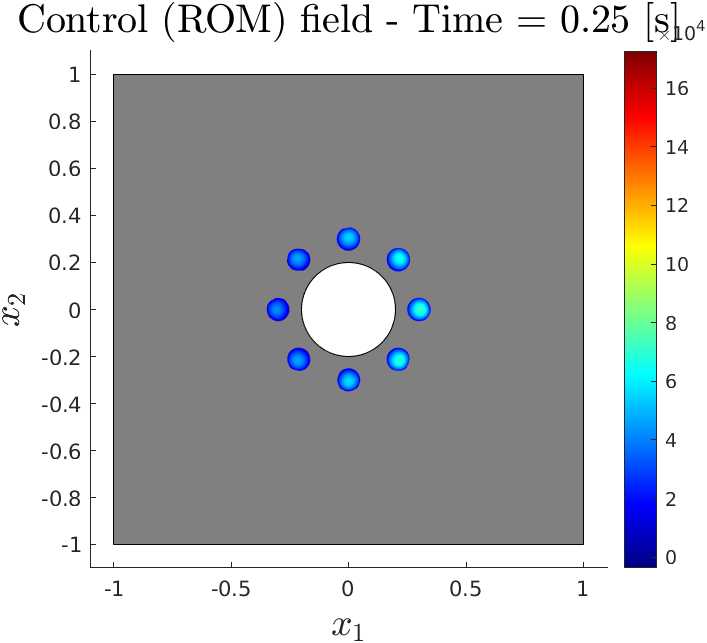}}
\subfigure{\includegraphics[width=0.3\textwidth]{images/circle_t/state_ROM_field_11.png}}
\subfigure{\includegraphics[width=0.3\textwidth]{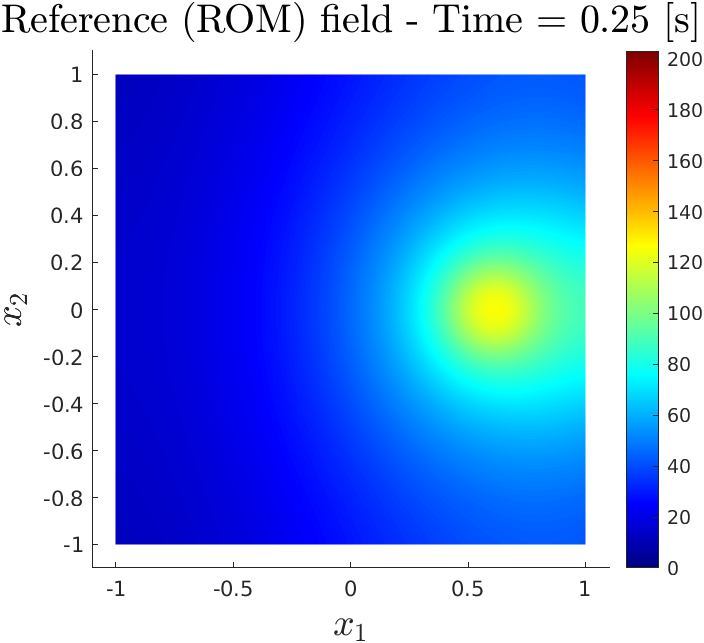}}
\caption{Reduction of the parametrized transient OCP with disconnected control domain, computational results for $\bs{\mu}_{t1} = [3.5 \, 10^4 \, 0] $ at time instance $t=0.25(\text{s})$.}
\label{circles_t_frame_1}
\end{figure}

\begin{figure}[h!]
\centering
\subfigure{\includegraphics[width=0.3\textwidth]{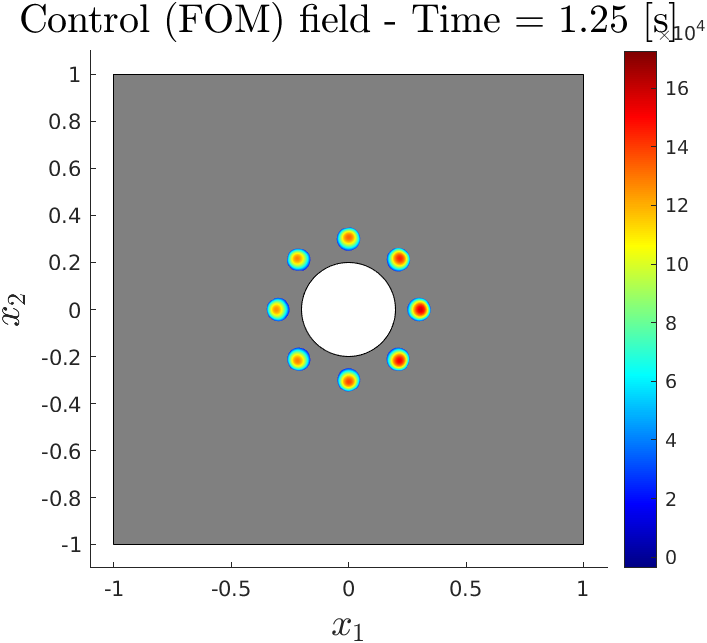}}
\subfigure{\includegraphics[width=0.3\textwidth]{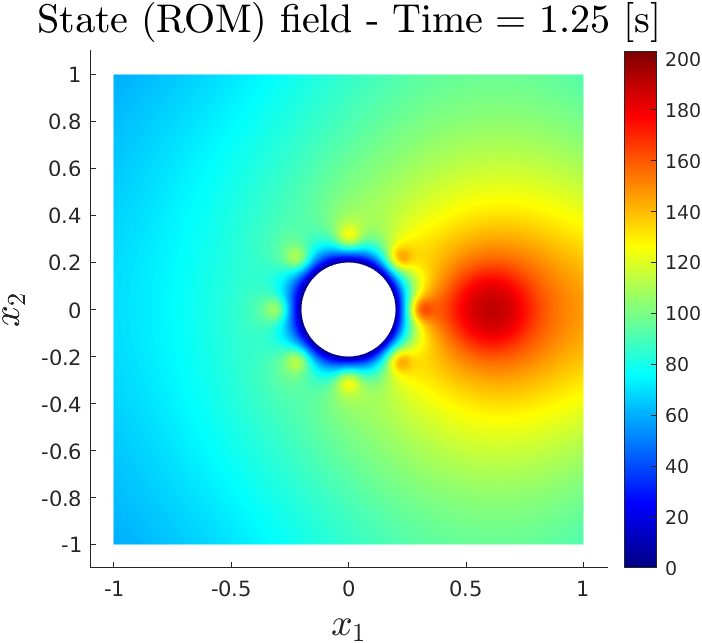}}
\subfigure{\includegraphics[width=0.3\textwidth]{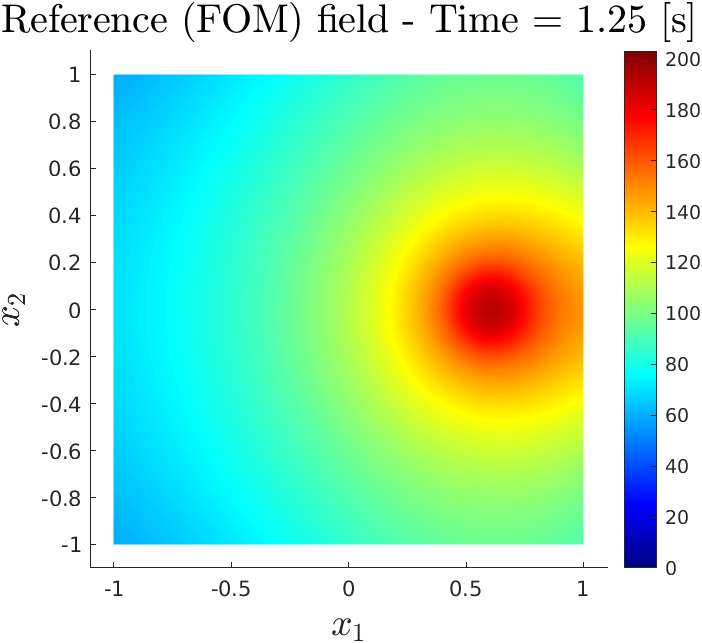}}
\subfigure{\includegraphics[width=0.3\textwidth]{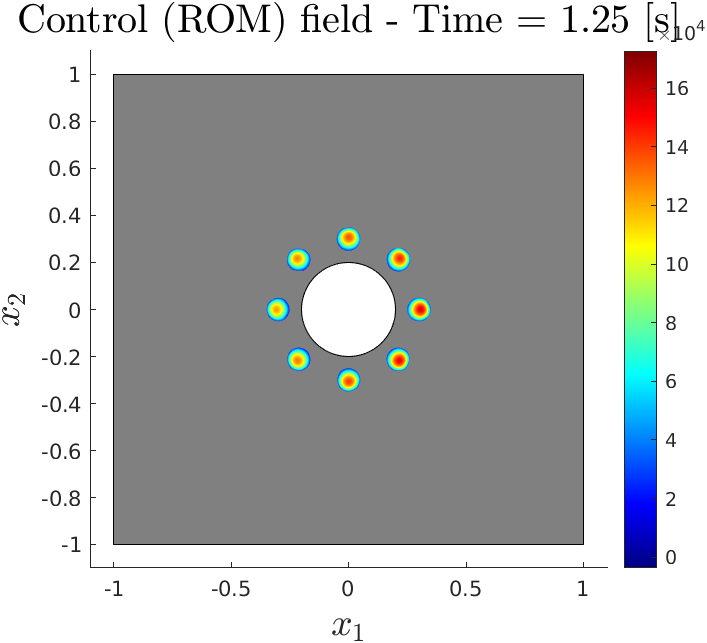}}
\subfigure{\includegraphics[width=0.3\textwidth]{images/circle_t/state_ROM_field_51.png}}
\subfigure{\includegraphics[width=0.3\textwidth]{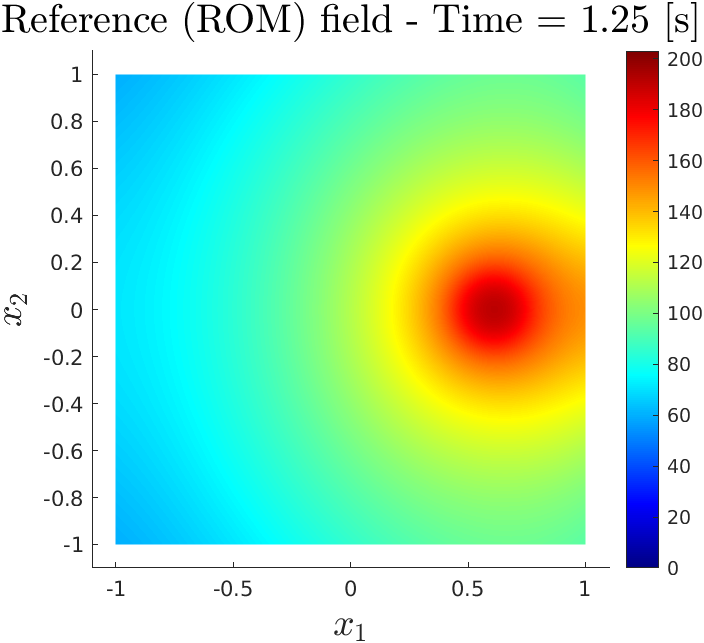}}
\caption{Reduction of the parametrized transient OCP with disconnected control domain, computational results for $\bs{\mu}_{t1} = [3.5 \, 10^4 \, 0] $ at time instance $t=1.25(\text{s})$.}
\label{circles_t_frame_2}
\end{figure}

\begin{figure}[h!]
\centering
\subfigure{\includegraphics[width=0.3\textwidth]{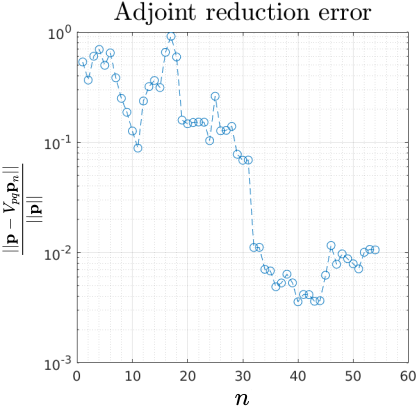}}
\subfigure{\includegraphics[width=0.3\textwidth]{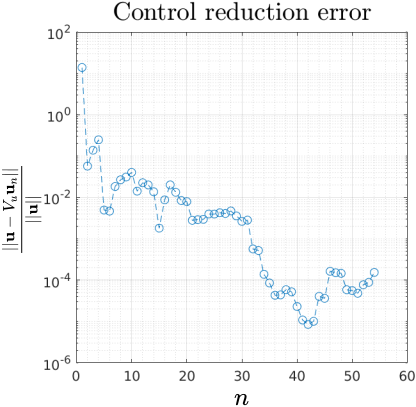}}
\subfigure{\includegraphics[width=0.3\textwidth]{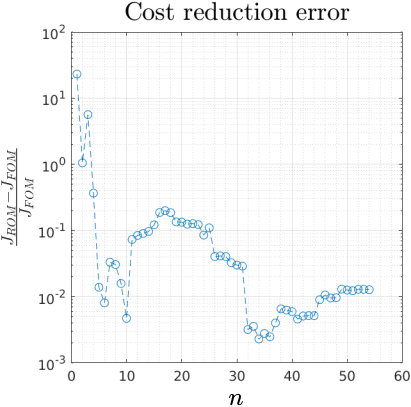}}
\subfigure{\includegraphics[width=0.3\textwidth]{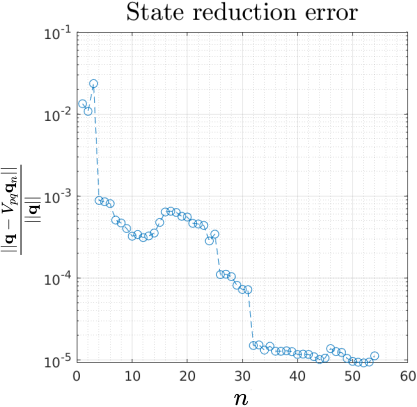}}
\subfigure{\includegraphics[width=0.3\textwidth]{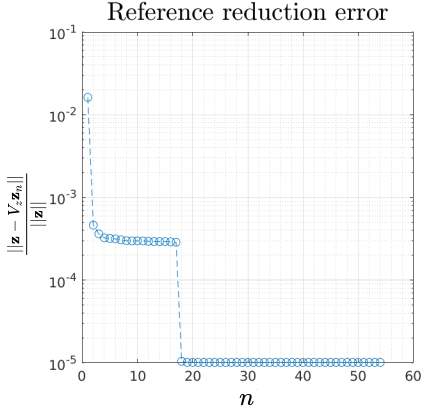}}
\subfigure{\includegraphics[width=0.3\textwidth]{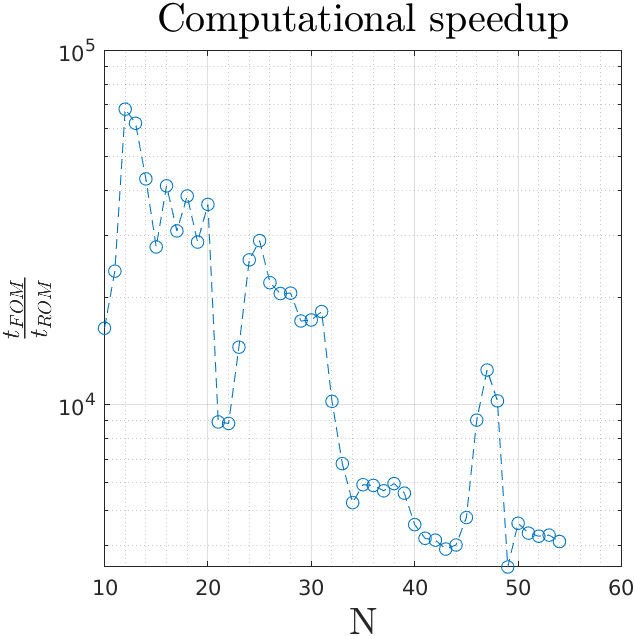}}
\caption{Relative errors between the ROM and FOM solutions for the transient problem of the disconnected control domain case. The ROM is able to achieve a computational speedup of $4000$ times without losing in accuracy of the optimal and reference solutions. The solution of the FOM takes 900 seconds while the ROM needs only 0.22 seconds.}
\label{circles_t_ROM}
\end{figure}

\noindent As for the previous test case, we first consider the reduction of the steady-state problem and then we tackle the transient one. As for the connected cloak, the reduction algorithm is able to achieve a huge computational speedup without any losing in accuracy. The same parameters range and number of snapshots is used to build the ROM and the reduction performances are shown in Figure \ref{circles_ss_ROM} and \ref{circles_t_ROM} for the steady-state and transient case respectively. The ROM and FOM field obtained for the disconnected case are shown in Figure \ref{circles_ss} for the steady-state case while three snapshots taken at time instances $t=0.25 (\text{s})$, $t=1.25 (\text{s})$  are shown in Figures \ref{circles_t_frame_1} and \ref{circles_t_frame_2}.

\begin{figure}[h!]
\centering
\subfigure{\includegraphics[width=0.3\textwidth]{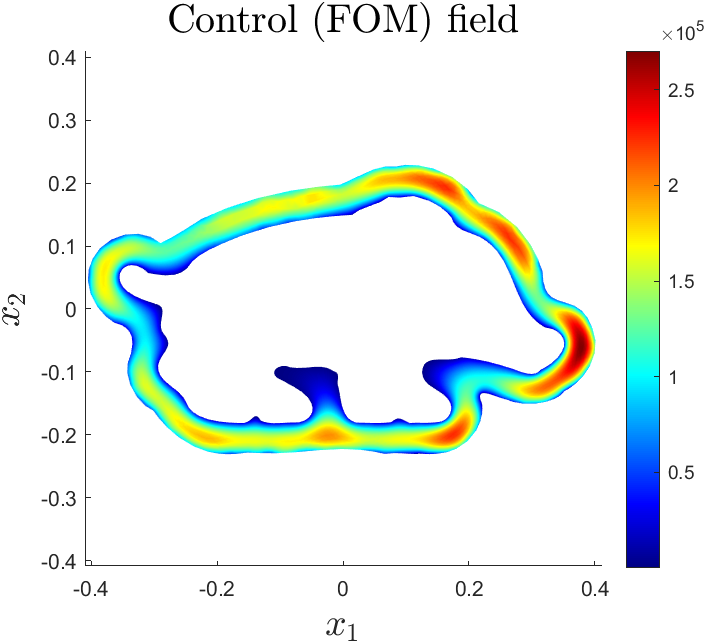}}
\subfigure{\includegraphics[width=0.3\textwidth]{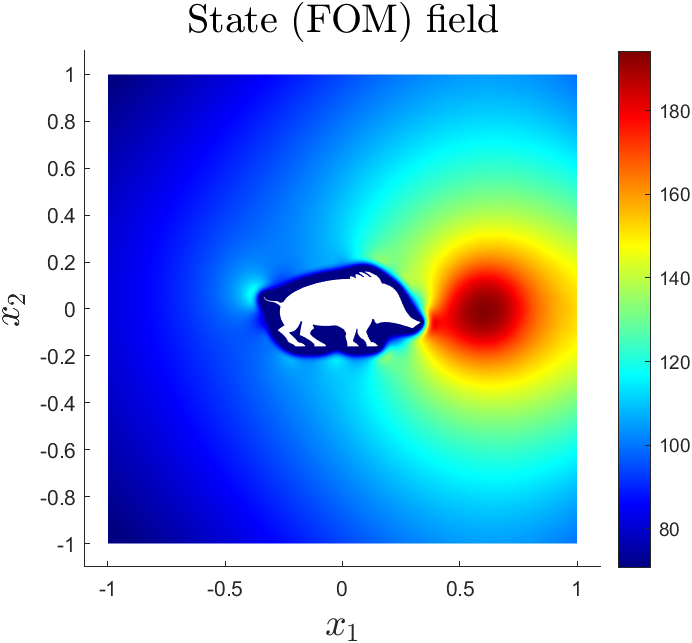}}
\subfigure{\includegraphics[width=0.3\textwidth]{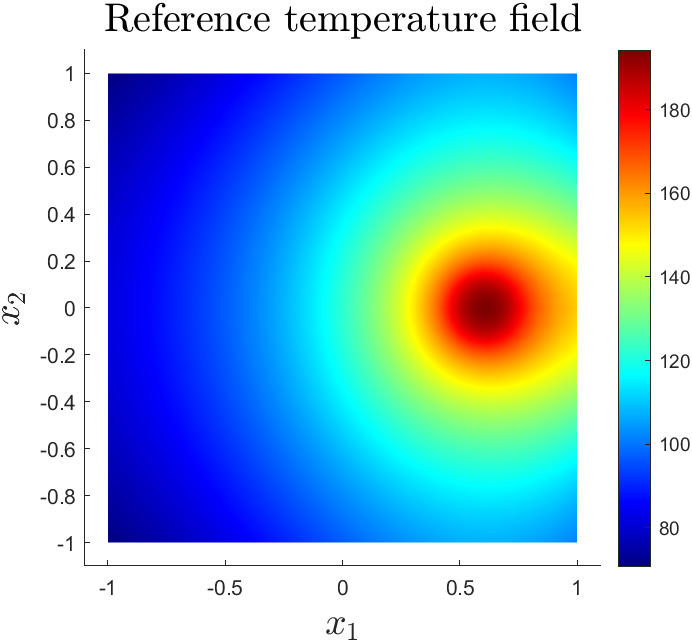}}
\subfigure{\includegraphics[width=0.3\textwidth]{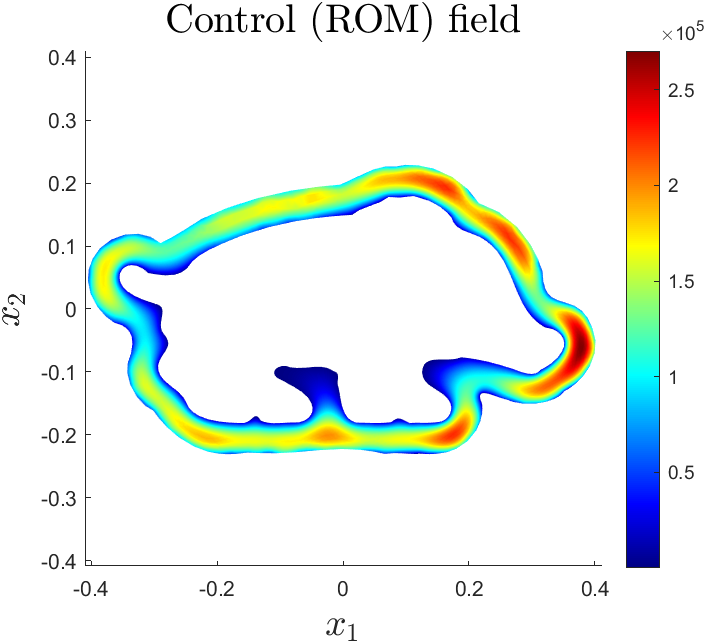}}
\subfigure{\includegraphics[width=0.3\textwidth]{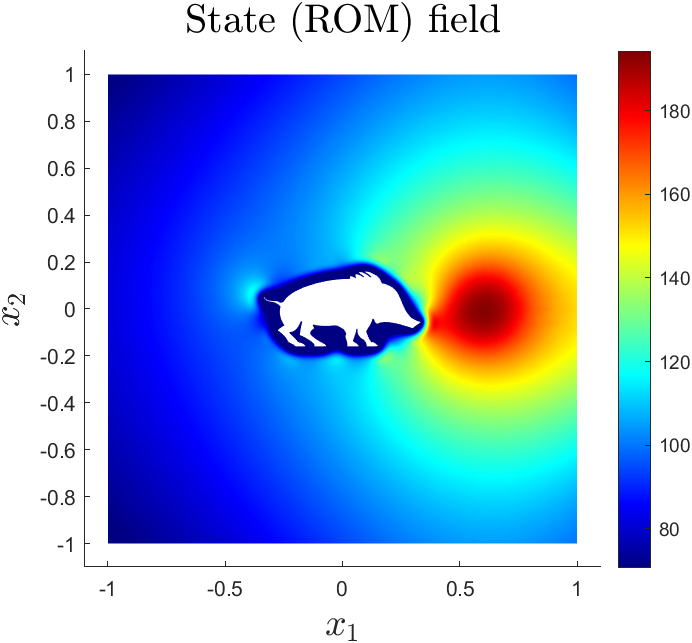}}
\subfigure{\includegraphics[width=0.3\textwidth]{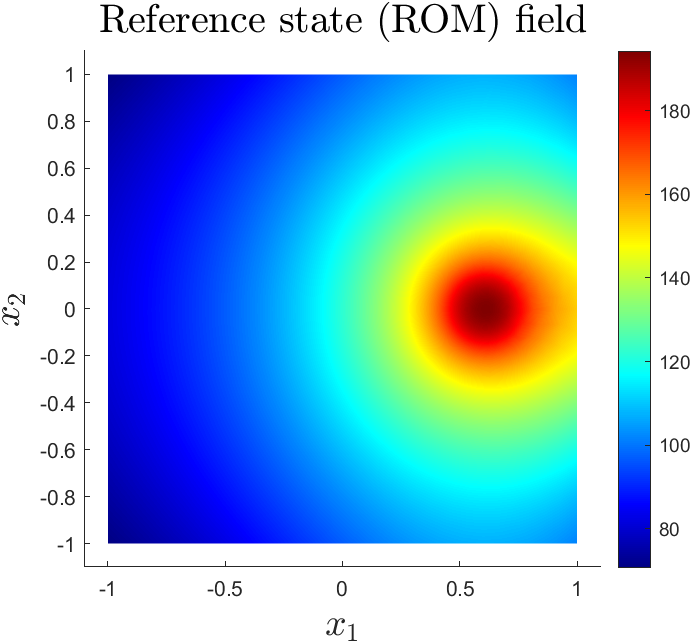}}
\caption{Reduction of the parametrized steady-state OCP concerning the complex shape case, computational results for $\bs{\mu}_{t1} = [\,3.5 \,\, 10^4 \,\, 0 \,] $.}
\label{scrofa_ss}
\end{figure}

\begin{figure}[h!]
\centering
\subfigure{\includegraphics[width=0.3\textwidth]{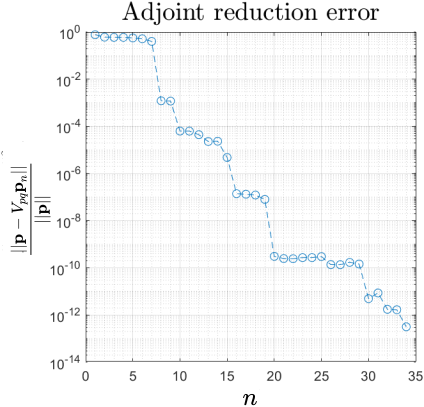}}
\subfigure{\includegraphics[width=0.3\textwidth]{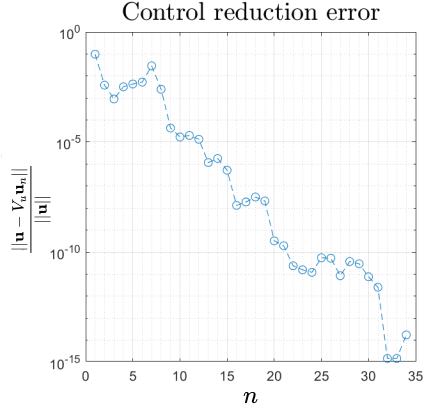}}
\subfigure{\includegraphics[width=0.3\textwidth]{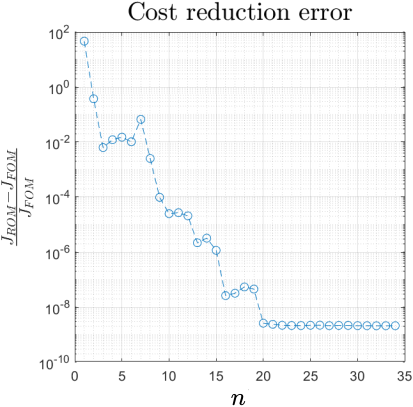}}
\subfigure{\includegraphics[width=0.3\textwidth]{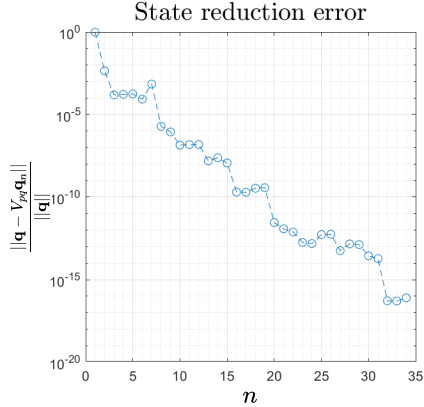}}
\subfigure{\includegraphics[width=0.3\textwidth]{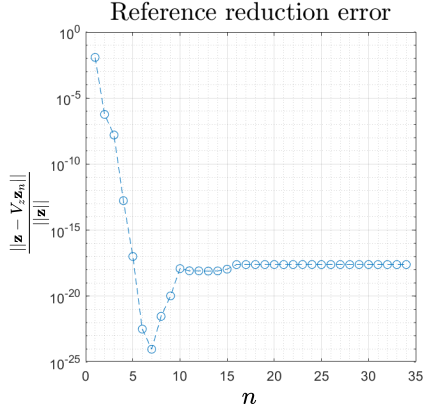}}
\subfigure{\includegraphics[width=0.3\textwidth]{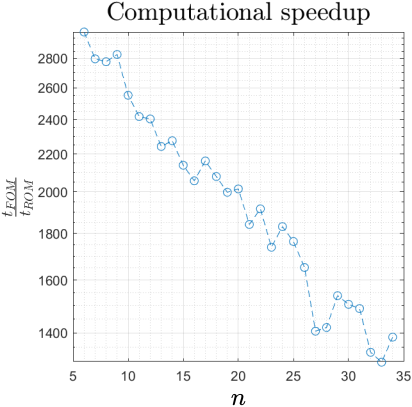}}
\caption{Relative errors between the ROM and FOM solutions for the complex shaped target. The ROM is able to achieve a computational speedup of $1385$ times without losing in accuracy of the optimal and reference solutions. Solution time for the projected linear system (\ref{big_ss}) is 0.78 milliseconds on a standard laptop computer while it takes 1.08 seconds to solve the FOM.}
\label{scrofa_ss_ROM}
\end{figure}

\begin{figure}[h!]
\centering
\subfigure{\includegraphics[width=0.3\textwidth]{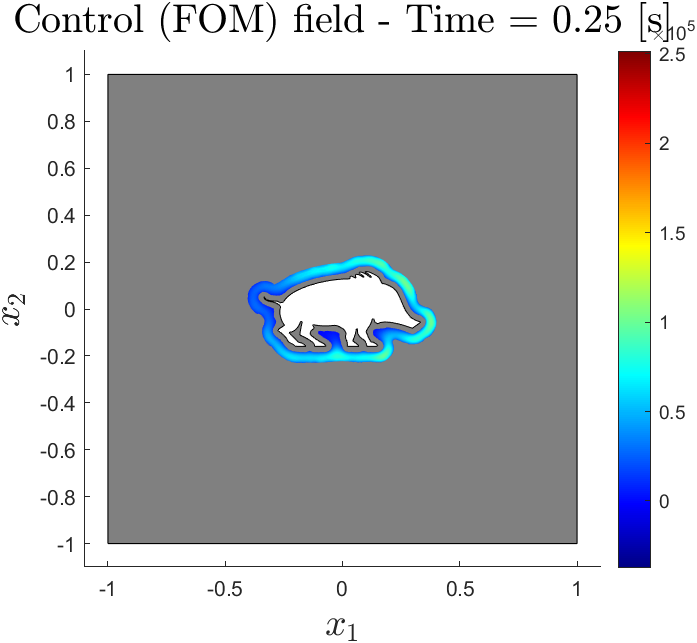}}
\subfigure{\includegraphics[width=0.3\textwidth]{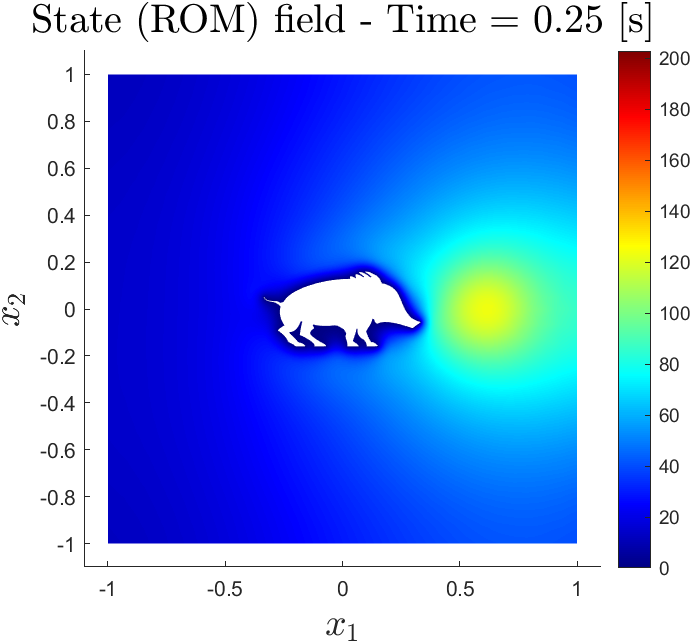}}
\subfigure{\includegraphics[width=0.3\textwidth]{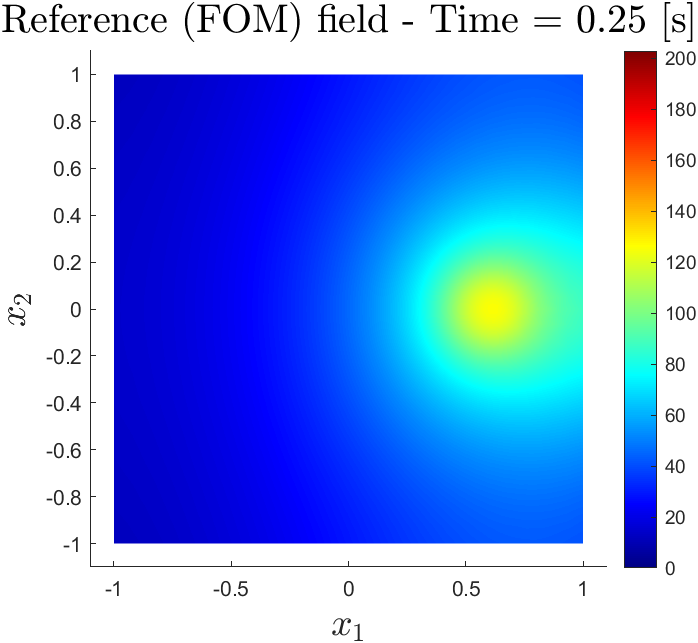}}
\subfigure{\includegraphics[width=0.3\textwidth]{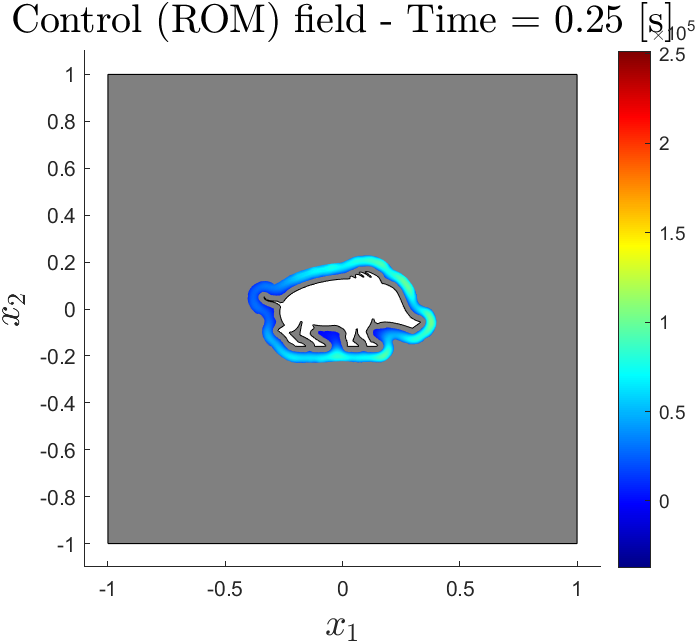}}
\subfigure{\includegraphics[width=0.3\textwidth]{images/scrofa_t/state_ROM_field_11.png}}
\subfigure{\includegraphics[width=0.3\textwidth]{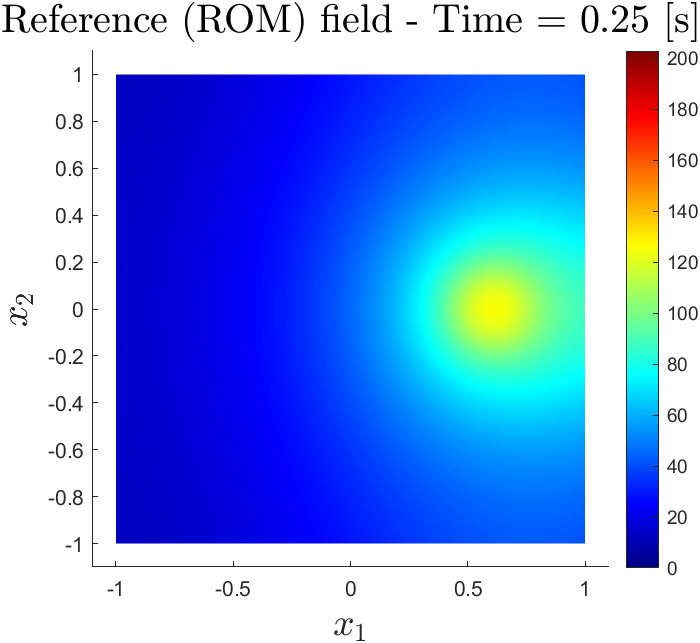}}
\caption{Reduction of the parametrized transient OCP with complex shaped obstacle and cloak, computational results for $\bs{\mu}_{t1} = [3.5 \, 10^4 \, 0] $ at time instance $t=0.25(\text{s})$.}
\label{scrofa_t_frame_1}
\end{figure}

\begin{figure}[h!]
\centering
\subfigure{\includegraphics[width=0.3\textwidth]{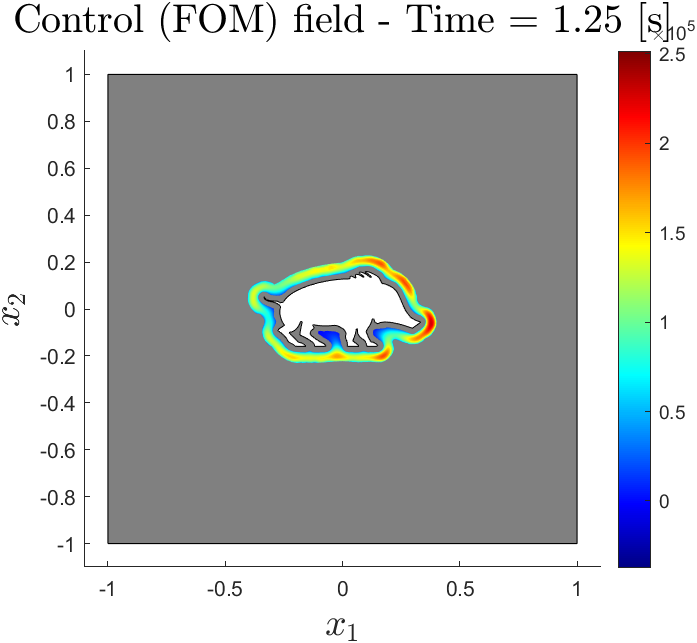}}
\subfigure{\includegraphics[width=0.3\textwidth]{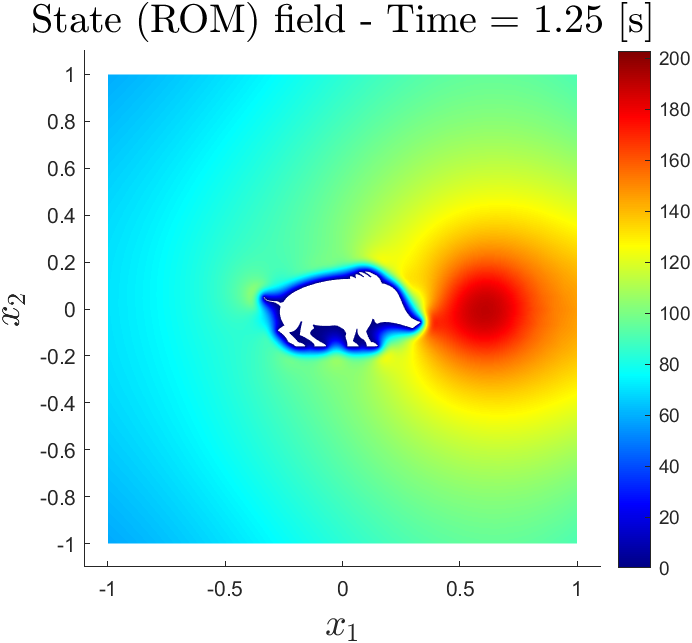}}
\subfigure{\includegraphics[width=0.3\textwidth]{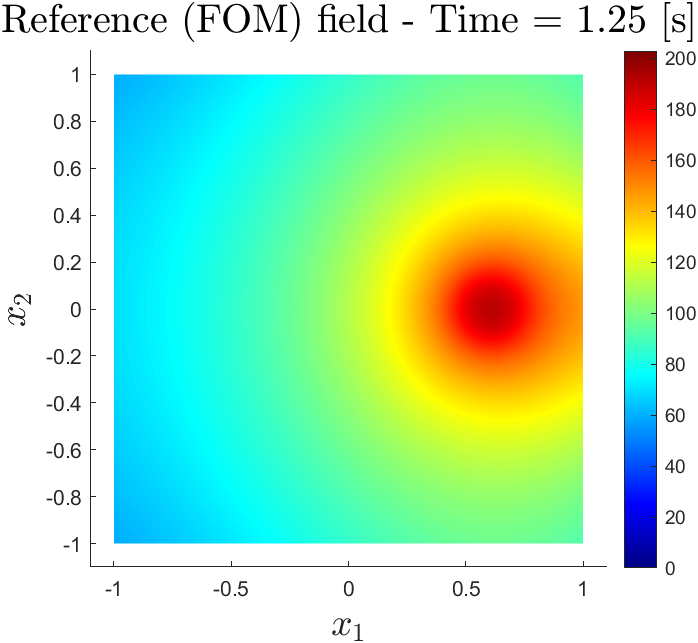}}
\subfigure{\includegraphics[width=0.3\textwidth]{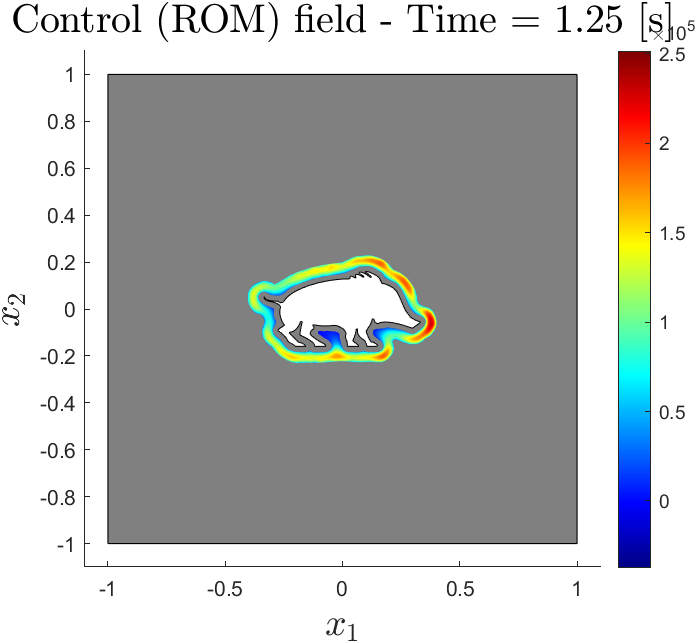}}
\subfigure{\includegraphics[width=0.3\textwidth]{images/scrofa_t/state_ROM_field_51.png}}
\subfigure{\includegraphics[width=0.3\textwidth]{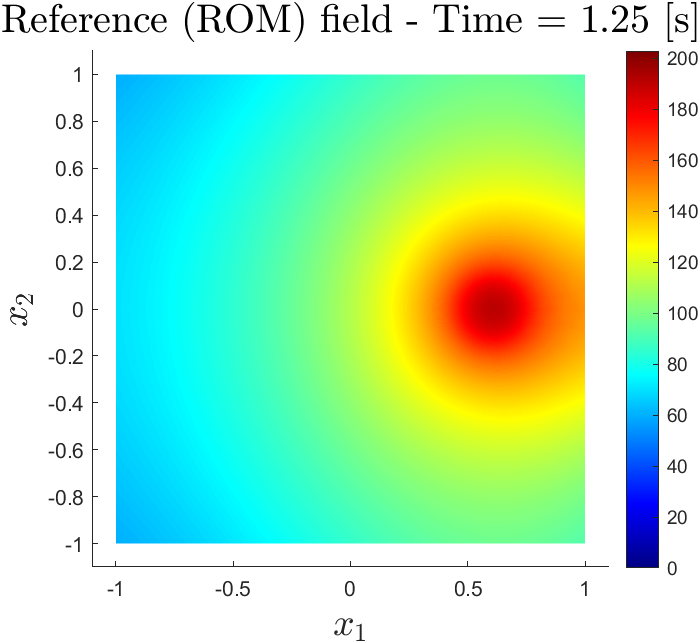}}
\caption{Reduction of the parametrized transient OCP with complex shaped obstacle and cloak, computational results for $\bs{\mu}_{t1} = [3.5 \, 10^4 \, 0] $ at time instance $t=1.25(\text{s})$.}
\label{scrofa_t_frame_2}
\end{figure}

\begin{figure}[h!]
\centering
\subfigure{\includegraphics[width=0.3\textwidth]{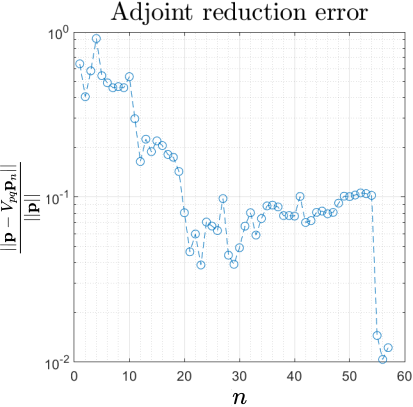}}
\subfigure{\includegraphics[width=0.3\textwidth]{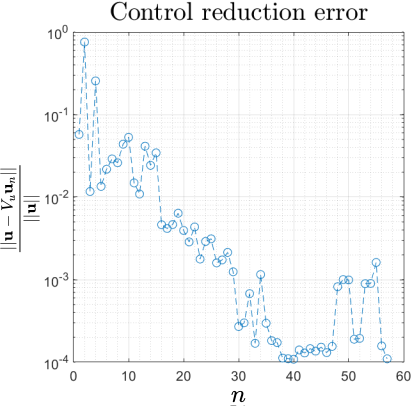}}
\subfigure{\includegraphics[width=0.3\textwidth]{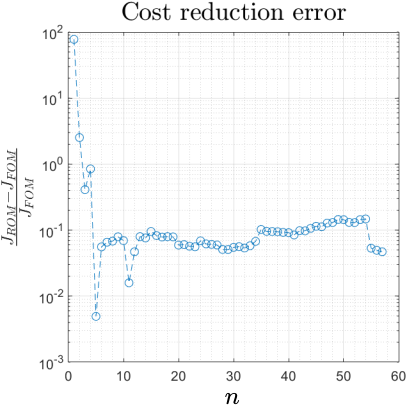}}
\subfigure{\includegraphics[width=0.3\textwidth]{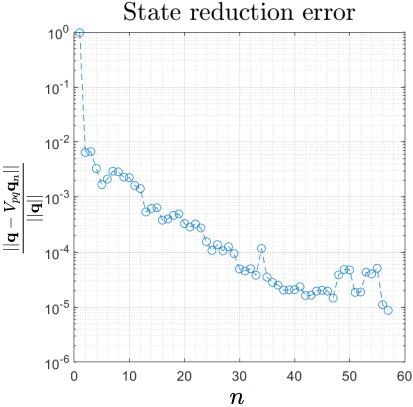}}
\subfigure{\includegraphics[width=0.3\textwidth]{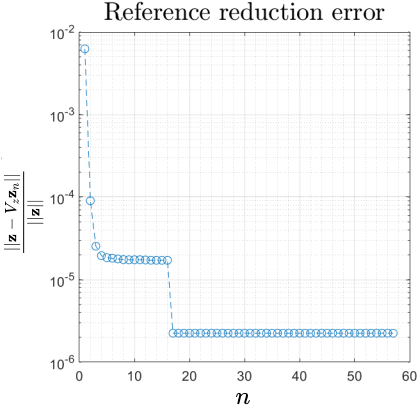}}
\subfigure{\includegraphics[width=0.3\textwidth]{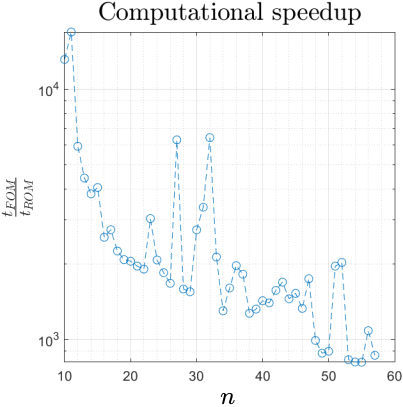}}
\caption{Relative errors between the ROM and FOM solutions for the transient problem of the complex shaped object. The ROM is able to achieve a computational speedup of $863$ times without losing in accuracy of the optimal and reference solutions. The solution of the FOM takes 394 seconds while the ROM needs only 0.45 seconds.}
\label{scrofa_t_ROM}
\end{figure}

As a final test case, we consider the thermal cloaking problem of a complex object where the cloaking sources are composed of a thin offset of the target, and thus almost replicate its complex shape. As a representative example, for the sake of dealing with complex shapes, we cloak the silhouette of the half-woollen boar, a legendary creature related to the foundation of the city of Milan. With the same ideas as for the previous test cases we consider the steady-state and transient regimes. Regarding the steady-state, the comparison between ROM and FOM is shown in Figure \ref{scrofa_ss}, while the convergence results are plotted in Figure \ref{scrofa_ss_ROM}. Due to the complexity of the cloak shape and of the target, the ROM solution time is roughly $260$ times faster with respect to the FOM one. However, the relative errors are still negligible. The transient reduction results are shown in Figure \ref{scrofa_t_frame_1} and \ref{scrofa_t_frame_2} where the snapshots are taken at the same time instances as for the previous cases. As shown in Figure \ref{scrofa_t_ROM}, the relative accuracy slightly worsens for the transient case but as it is evident from the figures that the reduction technique is extremely accurate. Note that due to the higher complexity of the problem, related with both the complex shape of the target and of the cloak, a richer basis is required to track the optimal variables and reproduce the OCP in the projected space. The computational speedup is 1385 with respect to the FOM solution. It takes roughly $1$ second to compute the solution of the FOM while the ROM needs only 0.78 milliseconds.

\begin{figure}[h!]
\centering
\subfigure[$\eta=0.999$]{\includegraphics[width=0.28\textwidth]{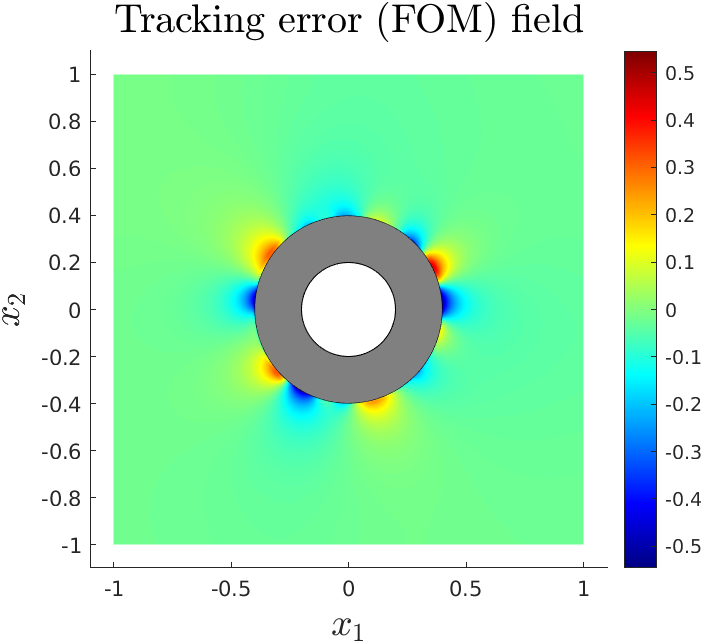}}
\subfigure[$\eta=0.989$]{\includegraphics[width=0.28\textwidth]{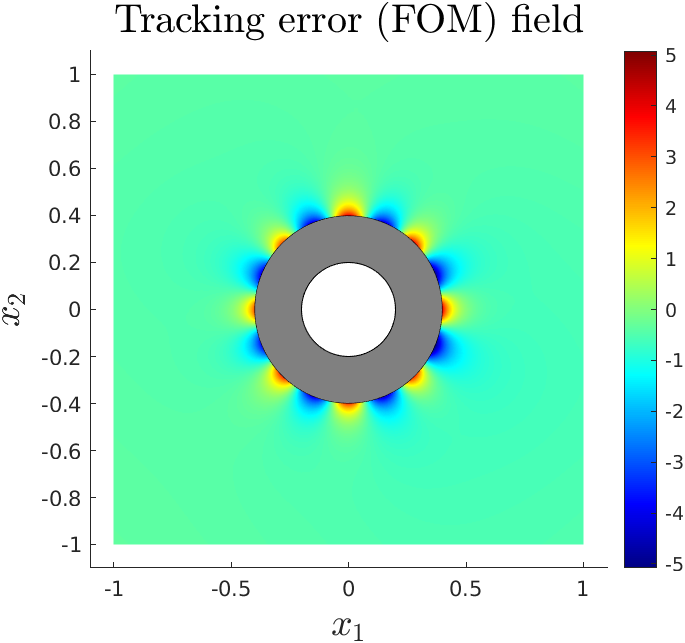}}
\subfigure[$\eta=0.966$]{\includegraphics[width=0.28\textwidth]{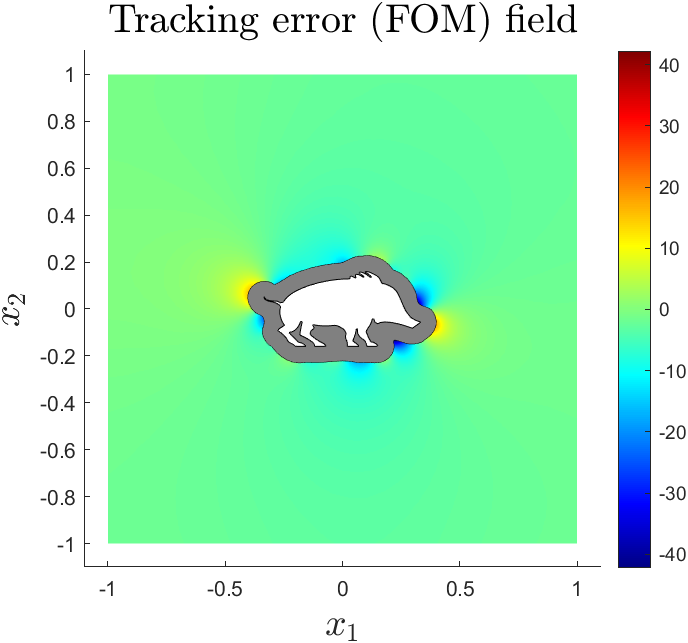}}
\caption{Tracking error fields of the Full Order Model for the layouts considered in this paper. The scale is symmetric about the maximum absolute value of the error. Water green color indicates relatively accurate tracking. The more complex is the problem the more the tracking error increases.}
\label{tracking_errors_ss}
\end{figure}

\noindent In Figure \ref{tracking_errors_ss}, we show the tracking error fields at steady-state for the three layouts considered in this work. The more complex the target to cloak and the shape or the cloak domain, the higher the tracking error. We measure the tracking performance computing the mean value in the $L^2(\Omega_{obs})$ sense of the squared difference between the optimal state and the reference state in the observation domain. We thus define the Mean Tracking Error (MTE) as:
\begin{equation*}
    MTE = \sqrt{\frac{\int_{\Omega_{obs}}  (q^{ss}(\x)-z^{ss}(\x))^2 \, d \Omega}{\int_{\Omega_{obs}} d \Omega}}.
\end{equation*}

\noindent The observation domain is defined as in Figure \ref{control_layout} for the circular and disconnected cloaks while it is the domain outside the cloak for the complex shape case. Note that the MTE is closely related to the first term in the cost functional of the OCP (\ref{ocp_formulation}), that is the term the optimization problem aims at minimizing while finding a weighted trade-off with the control norm. In order to compare the different cloaking layouts we compute the MTE for the uncontrolled case and define the cloaking efficiency as:
\begin{equation*}
\eta = \frac{|MTE-MTE^{\star}|}{MTE}
\end{equation*}
where $MTE^{\star}$ indicates the optimal variables while the absence of the superscript denotes the uncontrolled case. Note that a perfect cloak achieves $\eta=1$ while the uncontrolled case results in $\eta=0$. The resulting value of $\eta$ for each cloaking layout is shown in Figure \ref{tracking_errors_ss}.

\begin{figure}[h!]
\centering
\includegraphics[width=0.4\textwidth]{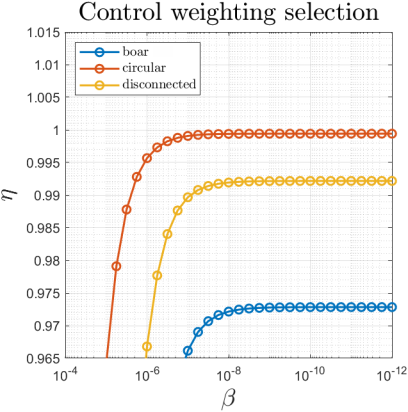}
\caption{Steady-state cloaking performances of the Full Order Model as a function of the control weighting parameter $\beta$ for the parameters set $\bs{\mu}_{t1} = [3.5 \, 10^4 \, 0]$. The gradient weighting parameter $\beta_g$ is set to $10^{-8} $. }
\label{beta_comparison}
\end{figure}

\noindent It is worthy to note that the cloaking error is higher in the proximity of the cloak and in particular where its profile is locally concave. Finally, we comment on the performance of the cloaking algorithm as the control weighting parameter varies since the optimization problem inherently depends on it. Figure \ref{beta_comparison} shows the steady-state performances for the three cases considered in this work as a function of the control parameter $\beta$. When the control weighting is sufficiently small, the performance of the cloak cannot be increased over a certain efficiency which depends implicitly on the complexity of the target to cloak and on the cloak shape.

\section{Conclusion}
\label{conclusion}
In this paper we have shown how to reformulate the design phase of the active thermal cloaking problem from an optimization standpoint. We recalled the theory of linear-quadratic OCP in the PDE settings and derived a system of necessary and sufficient conditions for optimality. Then, we developed a reduction framework to simultaneously treat the parametrized reference and optimal dynamics achieving a huge computational speedup without losing in accuracy, both in the transient and steady-state regimes and for a wide range of scenario parameters. The computational time needed to build the reduced basis offline is justified by the performance of the ROM that is able to compute the optimal solution for a time interval of $5$ seconds in approximately $0.2$ seconds. As a consequence, our reduction framework can be embedded in real-time applications where the reference field is measured or estimated. We leave as future work the matter of investigating such applications in, for example, a suitable Model Predictive Control (MPC) framework. Finally, we demonstrated that the optimal control framework allows to successfully treat complex and disconnected geometries both of the target to cloak and of the cloaking sources distribution.


\bibliographystyle{IEEEtran}
\bibliography{bibliography/biblist.bib}

\end{document}